\documentclass[11pt]{article}

\usepackage[utf8]{inputenc}

\usepackage{amssymb,xcolor,xspace}
\usepackage{amsmath,graphicx,comment,bm}
\usepackage{gastex}

\usepackage{hyperref}
\usepackage{array}

\usepackage{paralist}
\usepackage[linesnumbered,ruled,vlined]{algorithm2e}

\usepackage{geometry} 

\usepackage{fancyhdr} 
\usepackage{sectsty}
\allsectionsfont{\sffamily\mdseries\upshape} 

\def\cqfd{\skip10=\parfillskip\parfillskip=0pt
\enspace\hfill\symbolecqfd\par\parfillskip=\skip10\par\medskip}
\def\symbolecqfd{\rlap{$\sqcap$}$\sqcup$}

\newtheorem{theorem}{Theorem}[section]
\newtheorem{proposition}[theorem]{Proposition}
\newtheorem{lemma}[theorem]{Lemma}

\newtheorem{pro-fact}[theorem]{Fact}

\newtheorem{pro-example}[theorem]{Example}
\newenvironment{example}{\begin{pro-example}\rm}{\cqfd\end{pro-example}}

\newtheorem{pro-remark}[theorem]{Remark}
\newenvironment{remark}{\begin{pro-remark}\rm}{\cqfd\end{pro-remark}}

\newenvironment{preuve}{\rm \trivlist \item[\hskip \labelsep{\bf
Proof.}]}{\cqfd\endtrivlist}

\def\cqfd{\skip10=\parfillskip\parfillskip=0pt
\enspace\hfill\symbolecqfd\par\parfillskip=\skip10\par\medskip}
\def\symbolecqfd{\rlap{$\sqcap$}$\sqcup$}
\def\proof{\begin{preuve}}
\def\eop{\end{preuve}}

\let\phi\varphi

\DeclareMathOperator{\relab}{norm}

\def\inv{^{-1}}
\let\epsilon\varepsilon

\def \O {\mathcal{O}}

\def\E{\mathbb{E}}

\def\N{\mathbb{N}}

\def\PSL{\mathsf{PSL}_2(\Z)}

\def\Z{\mathbb{Z}}

\def\core{\mathsf{silh}}

\def\path{\textsf{path}}

\def\ella{\ell_2}
\def\ellb{\ell_3}
\def\ka{k_2}
\def\kb{k_3}

\def\rcrg{\mathtt{random\_cyclically\_reduced\_graph}}
\def\rsilh{\mathtt{random\_silhouette\_graph}}
\def\berna{\mathtt{bernoulli\_attempt}}
\def\shuffle{\mathsf{shuffle}}
\def\shift{\mathsf{shift}}

\title{Random generation of subgroups of the modular group with a fixed isomorphism type}

\author{
    Fr\'ed\'erique Bassino, \small{\url{bassino@lipn.fr}}\\
    \small{Universit\'e Sorbonne Paris Nord, LIPN, CNRS UMR 7030, F-93430 Villetaneuse, France}%
    \and
    Cyril Nicaud, \small{\url{cyril.nicaud@u-pem.fr}}\\
    \small{LIGM, Univ Gustave Eiffel, CNRS, ESIEE Paris, F-77454, Marne-la-Vallée, France}%
    \and
    Pascal Weil, \small{\url{pascal.weil@labri.fr}}\\
    \small{CNRS, ReLaX, IRL 2000, Siruseri, India} \\
    \small{Univ. Bordeaux, CNRS, Bordeaux INP, LaBRI, UMR 5800, F-33400 Talence, France}\thanks{%
    LaBRI, Univ. Bordeaux, 351 cours de la Lib\'eration, 33400 Talence, France.}
    }
    
\begin{document}

\maketitle

\begin{abstract}
We show how to efficiently count and generate uniformly at random finitely generated subgroups of the modular group $\PSL$ of a given isomorphism type. The method to achieve these results relies on a natural map of independent interest, which associates with any finitely generated subgroup of $\PSL$ a graph which we call its silhouette, and which can be interpreted as a conjugacy class of free finite index subgroups of $\PSL$. 
\end{abstract}




\section{Introduction}

The modular group $\PSL$ is a fundamental object in the field of modular forms and hyperbolic geometry. It is well-known that $\PSL$ is isomorphic to the free product of two cyclic groups, of order 2 and 3 respectively. That is,
$$\PSL = \langle a,b \mid a^2 = b^3 = 1\rangle.$$

The finitely generated subgroups of the modular group have been extensively studied and classified, leading to deep connections with various areas of mathematics, including number theory, algebraic geometry, and geometric group theory.  Much work has been devoted in particular to the combinatorial study of the \emph{finite index} subgroups of $\PSL$: exact enumeration results for the index $n$ subgroups (Dey, 1965 \cite{1965:Dey}; Stothers, 1978 \cite{1977:Stothers-Edinburgh}) and results on the asymptotic behavior of that number as $n$ tends to infinity (Newmann, 1976 \cite{1976:Newman}, M\"uller \& Schlage-Puchta, 2004 \cite{2004:MullerSchlage-Puchta} and others). Here, we deal instead with \emph{all} finitely generated subgroups of $\PSL$, without index restriction. 

The main purpose of this paper is to present
enumeration and random generation results for finitely generated subgroups of $\PSL$ of a given size and isomorphism type, where both measures are natural parameters.

Let us first make the notions underlying these results more explicit. Since $\PSL$ is the free product of a copy of $\Z_2$ and a copy of $\Z_3$, Kurosh's theorem (see, \textit{e.g.}, \cite{1956:Kurosh,1977:LyndonSchupp,1995:Rotman}) states that any finitely generated subgroup $H$ of $\PSL$ is isomorphic to a free product of $\ella$ copies of $\Z_2$, $\ellb$ copies of $\Z_3$ and $r$ copies of $\Z$: the \emph{isomorphism type} of $H$ is the triple $(\ella, \ellb, r)$. It is a natural parameter, which generalizes the rank in free groups.

Our results also refer to a notion of size for finitely generated subgroups of $\PSL$, that we now explain: each finitely generated subgroup $H$ of $\PSL$ can be represented uniquely by a finite edge-labeled graph $\Gamma(H)$, called its \emph{Stallings graph}.
Stallings graphs, and their effective construction, were first introduced by Stallings \cite{1983:Stallings} to represent finitely generated subgroups of free groups. The idea of using finite graphs to represent subgroups of infinite, non-free groups first appeared in work of Gersten and Short \cite{1991:GerstenShort,1991:Short}, Arzhantseva and Ol'shanski\u{\i} \cite{1996:ArzhantsevaOlshanskii,1998:Arzhantseva}, Gitik \cite{1996:Gitik} and Kapovich \cite{1996:Kapovich}. Markus-Epstein \cite{2007:Markus-Epstein} gave an explicit construction associating a graph with each subgroup of an amalgamated product of two finite groups, which is very close to the one used here. Here we follow the definition and construction of Kharlampovich, Miasnikov and Weil \cite{2017:KharlampovichMiasnikovWeil}. In a nutshell, the Stallings graph of a subgroup $H$ of $\PSL$ is the fragment of the Schreier (or coset) graph of $H$, spanned by the cycles at vertex $H$ reading a geodesic representative of an element of $H$, see Section~\ref{sec: stallings graph} for more details.

We take the number of vertices of $\Gamma(H)$ to be the \emph{size} of the subgroup $H$. In particular, there are only finitely many subgroups of a given size and we assume the uniform distribution on this finite set.

In \cite{2021:BassinoNicaudWeil} the authors counted the finitely generated subgroups of $\PSL$ by size and they showed how to generate uniformly at random a subgroup of a given size. They also computed the expected value of the isomorphism type of a random subgroup as a function of its size and proved a large deviations theorem for this isomorphism type. It follows that randomly generating a size $n$ subgroup of $\PSL$ will, with high probability, yield a subgroup whose isomorphism type is close to the average value. In particular, this algorithmic result does not help generate uniformly at random subgroups of a given size and isomorphism type.
The proof strategy to obtain these results was based on counting Stallings graphs and using the classical tools of analytic combinatorics \cite{2009:FlajoletSedgewick}, in particular the notion of exponential generating series.

In this paper, we use a completely different enumeration method for finitely generated subgroups of $\PSL$, to get a polynomial time random generation algorithm for subgroups of $\PSL$ of a given size and isomorphism type. It turns out that we can proceed with direct computations and we therefore avoid introducing generating series.  More precisely the proofs rely on a combination of graph decomposition techniques and combinatorial methods. As is classical in the field, these methods are used on labeled graphs (graphs equipped with a bijection from their vertex set to an initial segment of $\N$).

A key construction which occurs naturally in this approach is what we call the \emph{silhouetting} of the Stallings graph of a finitely generated subgroup of $\PSL$. It consists in a sequence of ``simplifications'' of the graph, leading (except in extremal cases) to a uniform degree loop-free graph, which represents a conjugacy class of finite index, free subgroups of $\PSL$.

The operation of silhouetting is not just useful for our enumeration and random generation purpose: it also has very interesting algebraic and probabilistic properties. As an example of the former, we establish that silhouetting preserves the free rank component of the isomorphism type of a subgroup (Proposition~\ref{prop: rank and small silhouette}). Probabilistic properties of the silhouetting operation, and their use in proving asymptotic properties of finitely generated subgroups of $\PSL$, are discussed in a separate paper \cite{2023:BassinoNicaudWeil-2}.

\paragraph{Organization of the paper}
Readers can find in Sections~\ref{sec: stallings graph} and~\ref{sec: combinatorial type} the precise definitions of the Stallings graph of a subgroup of $\PSL$ and its combinatorial type, and results from the literature relating this combinatorial information with algebraic properties of the subgroup such as its isomorphism type, its index or its freeness.

Section~\ref{sec: moves and silhouette} introduces combinatorial operations on Stallings graphs. Iterating these operations leads to so-called \emph{silhouette} graphs. The fine description of these operations is first exploited in Section~\ref{sec: counting} to give exact counting formulas for the number of subgroups of $\PSL$ of a given combinatorial or isomorphism type. 

In Section~\ref{sec: silhouette}, we show that the iteration of the operations defined in Section~\ref{sec: counting} is a confluent process (Proposition~\ref{prop: uniqueness silhouette}), which leads to defining the \emph{silhouette} of a given graph or subgroup. It is interesting to note that silhouetting preserves the free rank component of the isomorphism type of a subgroup (Proposition~\ref{prop: rank and small silhouette}).

Finally, Section~\ref{sec: random generation} uses the operations from Section~\ref{sec: moves and silhouette} in a different way to design an algorithm (which includes a rejection algorithm component) to efficiently generate uniformly at random a subgroup of a given size and isomorphism type.

\section{Preliminaries}\label{sec: rappels}

We work with the following presentation of the modular group:
$$\PSL = \langle a,b \mid a^2 = b^3 = 1\rangle.$$
The elements of $\PSL$ are represented by words over the alphabet $\{a,b,a\inv,b\inv\}$. Since $a\inv = a$ in $\PSL$, we can eliminate the letter $a\inv$ from this alphabet. Each element of $\PSL$ then has a unique shortest (or \emph{normal}, or \emph{geodesic}) representative, which is a freely reduced word without factors in $\{a^2, b^2, b^{-2}\}$. That is, the normal representatives are the words of length at most 1 and the words alternating letters $a$ and letters in $\{b,b\inv\}$.

\subsection{Stallings graph of a subgroup of $\PSL$}\label{sec: stallings graph}

The \emph{Schreier graph} (or \emph{coset graph}) of a subgroup $H$ of $\PSL$ is the graph whose vertices are the cosets $Hg$ of $H$ ($g\in \PSL$), with an $a$-labeled edge from $Hg$ to $Hga$ and a $b$-labeled edge from $Hg$ to $Hgb$, for every $g\in G$. We think of $b$-edges as $2$-way edges, reading $b$ in the forward direction and $b\inv$ in the backward direction. Since $a = a^2$, there is an $a$-edge from vertex $v$ to vertex $v'$ if and only there is one from $v'$ to $v$: as a result, we think of the $a$-edges as undirected edges, that can be traveled in either direction, each time reading $a$. A path in such an edge-labeled graph is called a \emph{cycle} if its initial vertex is equal to its final vertex. The sequence of edge labels along a path $p$ spells a word $w$ over alphabet $\{a,b,b\inv\}$, and we say that $w$ \emph{labels the path} $p$, or that $p$ \emph{reads} $w$.

Note that a word is in $H$ if and only if it labels a cycle at vertex $v_0 = H$ in the Schreier graph of $H$. The \emph{Stallings graph} of $H$, written $(\Gamma(H),v_0)$, is defined to be the fragment of the Schreier graph of $H$ spanned by the cycles at $v_0$ reading the geodesic representatives of the elements of $H$, rooted at $v_0$: that is, the subgraph of the Schreier graph of $H$ consisting of all the edges participating in a cycle at $v_0$ which reads a geodesic representation of an element of $H$, and of all the vertices adjacent to these edges. In particular, a word is in $H$ if and only if its geodesic representative labels a cycle in $\Gamma(H)$ at vertex $v_0$. We refer the reader to Remark~\ref{rk: historical} and to \cite{2017:KharlampovichMiasnikovWeil} for more details on these graphs. We note in particular that $H$ has a finite Stallings graph if and only if it is finitely generated, and that $\Gamma(H)$ is efficiently algorithmically computable if $H$ is given by a finite set of generators (words on the alphabet $\{a,b,b\inv\}$) \cite{2017:KharlampovichMiasnikovWeil,2021:BassinoNicaudWeil}.

\begin{example}\label{ex: Stallings graphs}
Figure~\ref{fig: stallings graphs} shows examples of Stallings graphs. To be more precise: the graphs in Figure~\ref{fig: stallings graphs} are \emph{labeled graphs}, meaning that their vertices are labeled by an initial segment of $\N$, see more details on this useful notion further down in Section~\ref{sec: labeled graphs}. The definition of Stallings graphs does not entail labeling vertices --- only designating a base vertex.
\begin{figure}[h!]
\centering
\begin{picture}(95,95)(0,-98)
\gasset{Nw=4,Nh=4}
\node(n0)(0.0,-0.0){$1$}

\node(n1)(6.0,-14.0){2}
\node(n2)(0.0,-28.0){3}
\node(n3)(34,-0.0){4}
\node(n4)(28,-14.0){5}
\node(n5)(34.0,-28.0){6}

\drawedge(n0,n1){$b$}

\drawedge(n1,n2){$b$}

\drawedge(n2,n0){$b$}

\drawedge[ELside=r](n3,n4){$b$}

\drawedge[ELside=r](n4,n5){$b$}

\drawedge[ELside=r](n5,n3){$b$}

\drawedge[AHnb=0](n0,n3){$a$}

\drawedge[AHnb=0](n1,n4){$a$}

\drawedge[AHnb=0](n2,n5){$a$}
\node(m1)(72.0,-4.0){$1$}
\node(m2)(78.0,-14.0){2}
\node(m3)(72.0,-24.0){3}
\node(m4)(52,-4.0){4}
\node(m5)(52,-24.0){5}
\node(m6)(98.0,-14.0){6}

\drawedge(m1,m2){$b$}

\drawedge(m2,m3){$b$}

\drawedge(m3,m1){$b$}

\drawedge[AHnb=0](m1,m4){$a$}

\drawedge[ELside=r](m4,m5){$b$}

\drawedge[AHnb=0](m5,m3){$a$}

\drawedge[AHnb=0](m2,m6){$a$}

\drawloop(m6){$b$}
\node(h1)(15.0,-40.0){$1$}
\node(h2)(45.0,-40.0){2}
\node(h3)(75.0,-40.0){3}
\node(h15)(105.0,-40.0){15}

\node(h4)(15.0,-64.0){4}
\node(h11)(15.0,-88.0){11}
\node(h16)(105.0,-64.0){16}
\node(h17)(105.0,-88.0){17}

\node(h12)(45.0,-88.0){12}
\node(h14)(75.0,-88.0){14}

\node(h10)(60.0,-56.0){10}
\node(h13)(30.0,-72.0){13}
\node(h18)(90.0,-72.0){18}

\node(h5)(4.0,-52.0){5}
\node(h6)(4.0,-76.0){6}
\node(h7)(-2.0,-86.0){7}
\node(h8)(4.0,-96.0){8}
\node(h9)(25.0,-96.0){9}

\drawedge(h2,h3){$b$}
\drawedge(h3,h10){$b$}
\drawedge(h10,h2){$b$}

\drawedge(h14,h17){$b$}
\drawedge[ELside=r](h17,h18){$b$}
\drawedge[ELside=r](h18,h14){$b$}

\drawedge(h11,h12){$b$}
\drawedge[ELside=r](h12,h13){$b$}
\drawedge[ELside=r](h13,h11){$b$}

\drawedge[ELside=r](h1,h4){$b$}
\drawedge(h4,h5){$b$}
\drawedge(h5,h1){$b$}

\drawedge(h6,h8){$b$}
\drawedge(h8,h7){$b$}
\drawedge(h7,h6){$b$}

\drawedge(h15,h16){$b$}

\drawedge[AHnb=0](h1,h2){$a$}
\drawedge[AHnb=0](h3,h15){$a$}
\drawedge[AHnb=0](h16,h17){$a$}
\drawedge[AHnb=0](h12,h14){$a$}
\drawedge[AHnb=0](h11,h4){$a$}
\drawedge[AHnb=0](h6,h5){$a$}
\drawedge[AHnb=0](h8,h9){$a$}
\drawedge[AHnb=0](h10,h13){$a$}

\drawloop[loopangle=140,AHnb=0](h7){$a$}
\drawloop[loopangle=0](h9){$b$}
\drawloop[loopangle=90,AHnb=0](h18){$a$}
\end{picture}
\caption{Top: the Stallings graphs of the subgroups
$H = \langle abab\inv, babab\rangle$ and $K = \langle abab, b^{ab\inv}\rangle$ of $\PSL$, where $g^h$ stands for $h\inv gh$. Bottom: the Stallings graph of 
  $L = \langle b^{ab\inv ab}, a^{bab}, a^{(ab)^3}, ab\inv abab\inv, (ab)^2(ab\inv)^3\rangle$.
In each case, the root is the vertex labeled 1.}\label{fig: stallings graphs}
\end{figure}
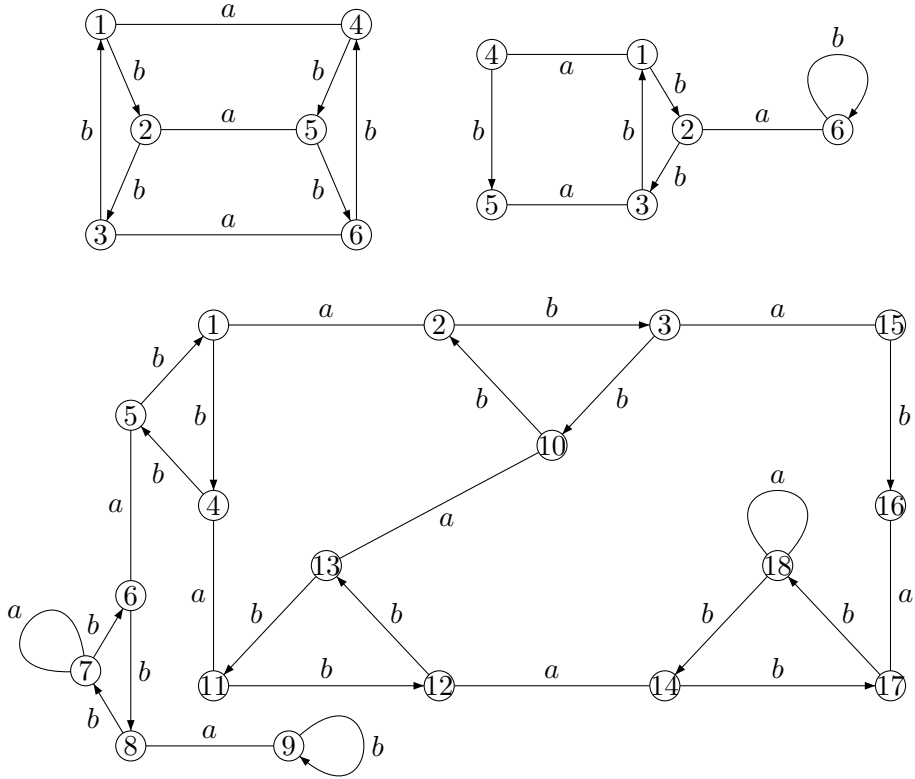
\end{example}

\begin{remark}\label{rk: historical}
  The definition of Stallings graphs given above is a generalization of that introduced by Stallings in 1983 \cite{1983:Stallings} for finitely generated subgroups of free groups, and a particular instance of the definition first introduced by Gitik \cite{1996:Gitik} in 1996 under the name of geodesic core, and systematized by Kharlampovich, Miasnikov and Weil \cite{2017:KharlampovichMiasnikovWeil} in 2017. Given a finitely presented group $G = \langle A\mid R\rangle$, a language $L$ of representatives for $G$ (a set of words over the alphabet $A \cup A\inv$) and a subgroup $H$, one considers the fragment of the Schreier graph of $H$ spanned by the $L$-representatives of the elements of $H$. It is effectively computable if $G$ is equipped with an automatic structure \cite{2017:KharlampovichMiasnikovWeil} and $H$ is $L$-quasi-convex. In the particular case
  where $G = \PSL$, we take $L$ to be the language of normal forms. It is well-known that $G$ is automatic with respect to this language, and that every finitely generated subgroup of $\PSL$ is quasi-convex. The algorithm to compute the Stallings graph of a subgroup, given a tuple of its generators, is quite straightforward, we refer the reader to \cite{2021:BassinoNicaudWeil} for an outline.
\end{remark}

It is immediate from the definition of Stallings graphs that $\Gamma(H)$ is connected and that its $a$-edges (respectively, $b$-edges) form a partial, injective map on the vertex set of the graph. Moreover, because $a^2 = b^3 = 1$, distinct $a$-edges are never adjacent to the same vertex: we distinguish therefore $a$-loops and so-called \emph{isolated $a$-edges}. Similarly, if we have two consecutive $b$-edges, say, from $v_1$ to $v_2$ and from $v_2$ to $v_3$, then $\Gamma(H)$ also has a $b$-edge from $v_3$ to $v_1$. Thus each $b$-edge is either a loop, or an \emph{isolated $b$-edge}, or a part of a $b$-triangle. Finally, every vertex except maybe the root vertex is adjacent to an $a$- and to a $b$-edge.

A rooted edge-labeled graph satisfying these conditions is called \emph{$\PSL$-reduced} and it is not difficult to see that every finite $\PSL$-reduced graph is the Stallings graph of a unique finitely generated subgroup of $\PSL$. That is, the mapping $H \mapsto (\Gamma(H),v_0)$ is a bijection between finitely generated subgroups of $\PSL$ and $\PSL$-reduced graphs.

An edge-labeled graph is said to be \emph{$\PSL$-cyclically reduced} if every vertex is adjacent to an $a$- and a $b$-edge or, equivalently, if it is $\PSL$-reduced when rooted at every one of its vertices. We also say that a finitely generated subgroup of $\PSL$ is $\PSL$-cyclically reduced if its Stallings graph is. 

\begin{example}\label{ex: small n}
The $\PSL$-cyclically reduced graphs $\Gamma$ with 1 or 2 vertices are represented in Figure~\ref{fig: 2-vertex cyclically reduced}. 
\begin{figure}[h!]
\centering
\begin{picture}(135,22)(3,-14)
\gasset{Nw=4,Nh=4}

\node(n0)(8.0,-5.0){}
\put(6,-15){$\Delta_1$}
\drawloop[loopangle=150,AHnb=0](n0){$a$}
\drawloop[loopangle=30](n0){$b$}

\node(n0)(33.0,-3.0){}
\node(n1)(53.0,-3.0){}
\put(41,-15){$\Delta_2$}
\drawedge[curvedepth=2](n0,n1){$b$}
\drawedge[curvedepth=2,AHnb=0](n1,n0){$a$}

\node(n0)(73.0,-5.0){}
\node(n1)(93.0,-5.0){}
\put(81,-15){$\Delta_3$}
\drawedge(n0,n1){$b$}
\drawloop[AHnb=0,loopangle=90](n0){$a$}
\drawloop[AHnb=0,loopangle=90](n1){$a$}

\node(n0)(113.0,-5.0){}
\node(n1)(133.0,-5.0){}
\put(121,-15){$\Delta_4$}
\drawedge[AHnb=0](n0,n1){$a$}
\drawloop[loopangle=90](n0){$b$}
\drawloop[loopangle=90](n1){$b$}
\end{picture}
\caption{\small All $\PSL$-cyclically reduced graphs with at most 2 vertices.}\label{fig: 2-vertex cyclically reduced}
\end{figure}
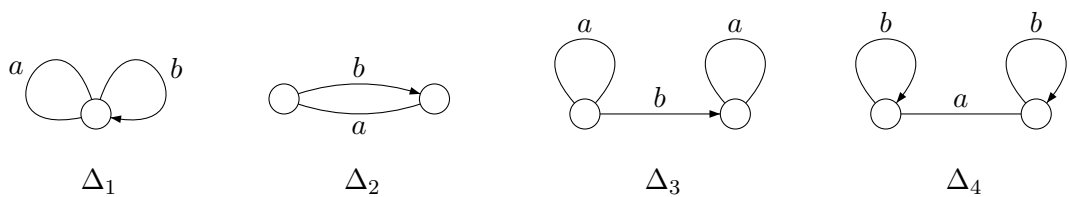

There is only one with 1 vertex, and three with 2 vertices. 
\end{example}

\subsection{Combinatorial type, isomorphism type of a subgroup of $\PSL$}\label{sec: combinatorial type}

The \emph{combinatorial type} of a $\PSL$-reduced graph $\Gamma$ is the tuple $(n,\ka,\kb,\ella,\ellb)$ where $n$ is the number of vertices of $\Gamma$, $\ka$ and $\kb$ are the numbers of isolated $a$- and $b$-edges, and $\ella$ and $\ellb$ are the numbers of $a$- and $b$-loops. We sometimes talk of the combinatorial type of a subgroup to mean the combinatorial type of its Stallings graph, and we refer to $n$ (the number of vertices of a $\PSL$-reduced graph) as the \emph{size} of the graph or even the \emph{size} of the subgroup. See \cite{2021:BassinoNicaudWeil} for a discussion of the possible combinatorial types.

Let us also record the following results (see, \emph{e.g.}, \cite[Lemma 2.3, Propositions 2.7, 2.9, 8.18 and Section 8.2]{2021:BassinoNicaudWeil}).

\begin{proposition}\label{prop: charact free and findex}
A subgroup $H\le \PSL$ has finite index if and only if its Stallings graph is $\PSL$-cyclically reduced and has combinatorial type of the form $(n,\ka,0,\ella,\ellb)$. It is free if and only if its combinatorial type is of the form $(n,\ka,\kb,0,0)$.

Free $\PSL$-cyclically reduced subgroups have even size. Free and finite index subgroups are $\PSL$-cyclically reduced and their size is a multiple of $6$. The combinatorial type of a free and finite index subgroup is of the form $(n, \frac12n, 0, 0, 0)$.

$\PSL$-cyclically reduced subgroups are conjugates if and only if the (unrooted) edge-labeled graphs underlying their Stallings graphs are isomorphic.
\end{proposition}

By Kurosh's classical theorem on subgroups of free groups (\textit{e.g.}, \cite[Proposition III.3.6]{1977:LyndonSchupp}), a subgroup $H$ of $\PSL$ is isomorphic to a free product of $r_2$ copies of $\Z_2$, $r_3$ copies of $\Z_3$ and a free group of rank $r$, for some non-negative integers $r_2,r_3,r$. The triple $(r_2,r_3,r)$, which characterizes $H$ up to isomorphism (but not up to an automorphism of $\PSL$) is called the \emph{isomorphism type} of $H$.
We record the following connection between the combinatorial and the isomorphism types of a subgroup \cite[Proposition 2.9]{2021:BassinoNicaudWeil}.

\begin{proposition}\label{prop: combinatorial vs. isomorphism}
Let $H$ be a subgroup of $\PSL$ of size at least $2$ and let $(n,\ka,\kb,\ella,\ellb)$ be the combinatorial type of $\Gamma(H)$.

If $\Gamma(H)$ is $\PSL$-cyclically reduced, the isomorphism type of $H$ is
$$\left(\ella,\ellb,1+\frac{n-2\kb-3\ella-4\ellb}6\right).$$

If $\Gamma(H)$ is not $\PSL$-cyclically reduced, the isomorphism type of $H$ is
\begin{align*}
\left(\ella,\ellb,\frac13+\frac{n-2\kb-3\ella-4\ellb}6\right) & \quad\text{if the base vertex is adjacent to an $a$-edge} \\
\left(\ella,\ellb,\frac12+\frac{n-2\kb-3\ella-4\ellb}6\right) & \quad\text{if the base vertex is adjacent to a $b$-edge.}
\end{align*}
\end{proposition}

\subsection{Labeled graphs}\label{sec: labeled graphs}

One of our objectives in this paper is to count subgroups by isomorphism type or by combinatorial type. Since subgroups are in bijection with $\PSL$-reduced graphs (their Stallings graphs), it is equivalent to count these graphs. For technical reasons, it is easier to count \emph{labeled graphs}, that is, graphs whose vertex set is equipped with a (labeling) bijection onto a set of the form $[n] = \{1,\dots,n\}$\footnote{It is important to distinguish this notion of labeling, which injectively assigns an integer to each vertex, from the edge labeling used so far, where each edge is labeled by either the order 2 generator $a$ of $\PSL$, or by its order 3 generator $b$ and each path is labeled by a word.}. The graphs in Figure~\ref{fig: stallings graphs} are in fact labeled graphs.

\begin{example}\label{ex: labelings of small graphs}
The $\PSL$-cyclically reduced graphs $\Delta_2$ and $\Delta_3$ in Example~\ref{ex: small n} admit two distinct labelings, while $\Delta_1$ and $\Delta_4$ have only one. We record here what we call the \emph{preferred} labeling of $\Delta_2$, where the $b$-edge goes from 1 to 2, see Figure~\ref{fig: canonical Delta2}.
\begin{figure}[htbp]
\centering
\begin{picture}(35,9)(3,-6)
\gasset{Nw=4,Nh=4}
\node(n0)(8.0,-3.0){1}
\node(n1)(33.0,-3.0){2}
\drawedge[curvedepth=2](n0,n1){$b$}
\drawedge[curvedepth=2,AHnb=0](n1,n0){$a$}
\end{picture}
\caption{\small The preferred labeling of $\Delta_2$}\label{fig: canonical Delta2}
\end{figure}
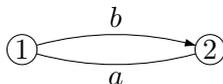
\end{example}

Since we are going to count graphs, rooted or not, labeled or not, it is important to clarify that we consider these combinatorial objects up to isomorphism. Concretely, if $\Gamma$ and $\Gamma'$ are graphs, an isomorphism from $\Gamma$ to $\Gamma'$ is a pair of bijections $\phi = (\phi_V, \phi_E)$ from the vertex set of $\Gamma$ to the vertex set of $\Gamma'$ and from the edge set of $\Gamma$ to the edge set of $\Gamma'$, respectively, which preserve the incidence relation (that is, if $\Gamma$ has an edge $e$ from vertex $v$ to vertex $w$, then $\phi_E(e)$ is an edge from $\phi_V(v)$ to $\phi_V(w)$. If $\Gamma$ and $\Gamma'$ are rooted, then $\phi_V$ must also map the root of $\Gamma$ to the root of $\Gamma'$. If, finally, $\Gamma$ and $\Gamma'$ are edge-labeled, then $\phi_E$ must also preserve these labels.

It is important to note that, an $n$-vertex $\PSL$-reduced graph admits exactly $n!$ distinct labeling functions. Let indeed $v_0$ be the root of $\Gamma$ and let us fix a total order on the alphabet $\{a,b,b\inv\}$. Assigning to each vertex $v$ the lexicographically least geodesic word labeling a path from $v_0$ to $v$, yields a total order on the vertex set of $\Gamma$. A labeling of $\Gamma$ is therefore equivalent to a permutation of $[n]$. Another way of formulating this observation is that a $\PSL$-reduced graph admits no non-trivial automorphism.

\section{The silhouetting operation on $\PSL$-cyclically reduced graphs}\label{sec: moves and silhouette}

We will see in Sections~\ref{sec: reduction to cyclically reduced} and~\ref{sec: random generation} that counting and randomly generating subgroups of $\PSL$ reduces to counting and randomly generating labeled $\PSL$-cyclically reduced graphs. 
Before we embark on this task, we introduce a combinatorial construction on this class of graphs.

More precisely, we define in Section~\ref{sec: moves} certain moves on a labeled $\PSL$-cyclically reduced graph, depending on its geometry. They are used in Section~\ref{sec: couting cyclically reduced} to count subgroups of $\PSL$ and in Section~\ref{sec: random generation} to randomly generate them. They also bring to the fore an interesting structure associated with a $\PSL$-cyclically reduced graph, which we call its \emph{silhouette}. Some of its algebraic and combinatorial properties are discussed in Section~\ref{sec: silhouette}.

Very roughly speaking, these moves ``simplify'' a labeled $\PSL$-cyclically reduced graph by first iteratively removing all loops and the paths that lead to them, until we are left (except in degenerate cases) with a graph which consists only of $b$-triangles and paths connecting them. The process then ``simplifies'' these connecting paths so that the resulting graph consists only of $b$-triangles connected by isolated $a$-edges. As we know, such graphs represent conjugacy classes of free finite index subgroups (see Proposition~\ref{prop: charact free and findex}).

For technical reasons, we use the notion of \emph{weakly labeled} graphs \cite[Definition II.1]{2009:FlajoletSedgewick}: if $\Gamma$ is a $\PSL$-cyclically reduced graph of size $m$, a \emph{weak labeling} of $\Gamma$ is an injective map from the vertex set of $\Gamma$ to $[n]$, where $n$ is an integer at least equal to $m$. To lighten up notation, we often abusively identify the vertices of a weakly labeled graph with their labels. We also abusively write $\Delta_i$ ($i = 1,2,3,4$) for any weakly labeled version of the graphs in Example~\ref{ex: small n}.

Observe that a weak labeling $\alpha$ of $\Gamma$ gives rise to a labeling of $\Gamma$ by a uniquely defined order-preserving bijection from the range of $\alpha$ to $[m]$. The labeled graph obtained this way is called the \emph{normalization}\footnote{This operation is called \emph{reduction} in \cite[Section II.2.1]{2009:FlajoletSedgewick}.} of $\Gamma$, denoted by $\relab(\Gamma)$.

\subsection{Moves on a labeled $\PSL$-cyclically reduced graph}\label{sec: moves}

Here we define so-called $\lambda_3$-, $\lambda_{2,1}$-, $\lambda_{2,2}$-, $\kappa_3$- and \emph{exceptional} moves on weakly labeled $\PSL$-cyclically reduced graphs\footnote{The denomination of $\lambda_3$-move is chosen because these moves deal with loops labeled by the order 3 generator $b$, which are counted by the parameter $\ellb$. Similar justifications hold for the moves that deal with $a$-loops (counted by $\ella$) and with isolated $b$-edges (counted by $\kb$).}.  But for the exceptional moves, each of these moves deletes vertices from the input graph without changing their label, so that the resulting graph is, again, weakly labeled.

In the following, $\Gamma$ is a weakly labeled $\PSL$-cyclically reduced graph with combinatorial type $\bm\tau = (n,\ka,\kb,\ella,\ellb)$.

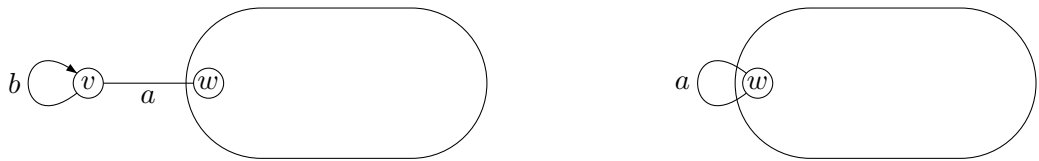
\begin{figure}[htbp]
\centering
\begin{picture}(130,22)(-2,-22)
\gasset{Nw=4,Nh=4}
\drawoval(37,-12,40,20,12)
\node(n0)(20.0,-12.0){$w$}
\node(n1)(4.0,-12.0){$v$}
\drawedge[AHnb=0](n0,n1){$a$}
\drawloop[loopangle=180,loopdiam=6](n1){$b$}

\drawoval(110,-12,40,20,12)
\node(n2)(93.0,-12.0){$w$}
\drawloop[loopangle=180,loopdiam=6,AHnb=0](n2){$a$}
\end{picture}
\caption{\small $(\lambda_3,v,w)$-move%
}\label{fig: tb differenciation}
\end{figure}

\paragraph{$\lambda_3$-moves}
If $\Gamma$ has a $b$-loop at vertex $v$ (in fact, at the vertex labeled $v$) and there is an isolated $a$-edge between $v$ and a distinct vertex $w$, then the \emph{$(\lambda_3,v,w)$-move} consists in deleting vertex $v$ and the adjacent edges, and adding an $a$-loop at vertex $w$. The resulting weakly labeled graph $\Delta$ (see Figure~\ref{fig: tb differenciation}) is $\PSL$-cyclically reduced and has combinatorial type $\bm\tau + \bm\lambda_3$, where $\bm\lambda_3 =  (-1,-1,0,1,-1)$.

\begin{lemma}\label{lm: claim 1}
Suppose that $n \ge 2$ and $\ellb > 0$. The $\lambda_3$-moves establish a bijection from the set of pairs $(\Gamma,\ell)$ with $\Gamma$ a labeled $\PSL$-cyclically reduced graph of combinatorial type $\bm\tau$ and $\ell$ a $b$-loop in $\Gamma$, to the set of triples $(\Delta,\ell',v)$ formed by a labeled $\PSL$-cyclically reduced graph $\Delta$ with combinatorial type $\bm\tau + \bm\lambda_3$, an $a$-loop $\ell'$ in $\Delta$ and an integer $v\in [n]$. 
\end{lemma}

\proof
Given a pair $(\Gamma,\ell)$ with $\Gamma$ a labeled $\PSL$-cyclically reduced graph of combinatorial type $\bm\tau$ and $\ell$ a $b$-loop in $\Gamma$, we associate to it the triple $(\Delta,v,w)$ --- and we write $(\Gamma,\ell) \longmapsto (\Delta,\ell',v)$ --- defined as follows: $v$ is the vertex carrying the loop $\ell$ in $\Gamma$; since $n > 1$, $v$ is adjacent to an isolated $a$-edge and we let $w$ be the other end of that $a$-edge; finally, we let $\Delta'$ be the weakly labeled graph obtained from $\Gamma$ by a $(\lambda_3,v,w)$-move, $\Delta = \relab(\Delta')$ and $\ell'$ be the $a$-loop in $\Delta$ at the vertex labeled $w$ in $\Delta'$.

Conversely, let $\Delta$ be a labeled $\PSL$-cyclically reduced graph with combinatorial type $\bm\tau + \bm\lambda_3$, let $\ell'$ be an $a$-loop in $\Delta$ and let $v \in [n]$. Let $\Delta'$ be the weakly labeled graph obtained from $\Delta$ by ``making space for $v$'', that is, by incrementing the labels of all the vertices greater than or equal to $v$. Finally, let $w$ be the label of the vertex of $\Delta'$ carrying the loop $\ell'$. Now let $\Gamma$ be the graph obtained from $\Delta'$ by deleting the loop $\ell'$ and adding vertex $v$, an isolated $a$-edge between $v$ and $w$ and a $b$-loop at $v$: it is directly verified that $\Gamma$ is properly labeled, of combinatorial type $\bm\tau$, and that $(\Gamma,\ell) \longmapsto (\Delta,\ell',v)$.
\eop

\paragraph{$\lambda_2$-moves}
Let $\Gamma$ be a weakly labeled $\PSL$-cyclically reduced graphs of size $n \ge 3$ and let $v$ be a vertex carrying an $a$-loop. Two situations occur, depending on whether $v$ sits on a $b$-triangle or not, giving rise to two flavors of $\lambda_2$-moves.

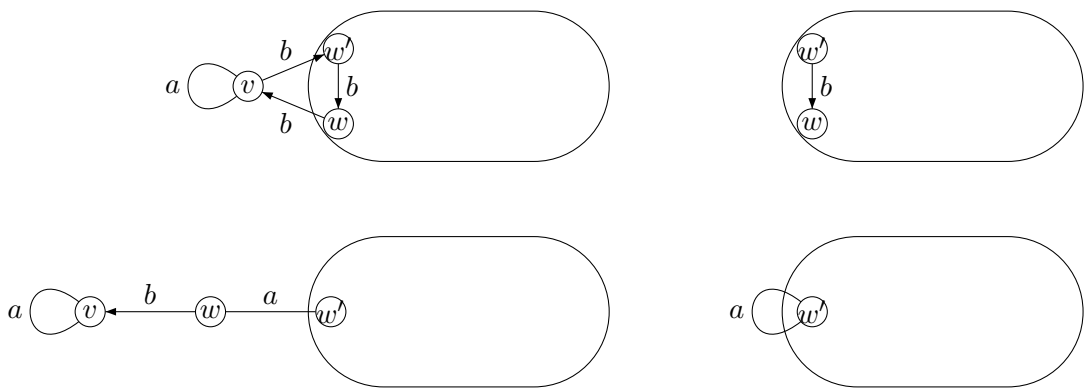
\begin{figure}[htbp]
\centering
\begin{picture}(130,52)(-2,-52)
\gasset{Nw=4,Nh=4}
\drawoval(49,-12,40,20,12)
\node(n0)(33.0,-7.0){$w'$}
\node(n00)(33.0,-17.0){$w$}
\node(n1)(21.0,-12.0){$v$}
\drawedge(n00,n1){$b$}
\drawedge(n1,n0){$b$}
\drawedge(n0,n00){$b$}
\drawloop[loopangle=180,loopdiam=6,AHnb=0](n1){$a$}

\drawoval(112,-12,40,20,12)
\node(n2)(96.0,-7.0){$w'$}
\node(n22)(96.0,-17.0){$w$}
\drawedge(n2,n22){$b$}

\drawoval(49,-42,40,20,12)
\node(n3)(32.0,-42.0){$w'$}
\node(n4)(16.0,-42.0){$w$}
\node(n5)(0.0,-42.0){$v$}
\drawedge[ELside=r](n4,n5){$b$}
\drawedge[ELside=r,AHnb=0](n3,n4){$a$}
\drawloop[loopangle=180,loopdiam=6,AHnb=0](n5){$a$}

\drawoval(112,-42,40,20,12)
\node(n6)(96.0,-42.0){$w'$}
\drawloop[loopangle=180,loopdiam=6,AHnb=0](n6){$a$}
\end{picture}
\caption{\small Above: $(\lambda_{2,1},v,w')$-move. Below: $(\lambda_{2,2}, v\leftarrow w,w')$-move}\label{fig: ta differenciation}
\end{figure}

If $v$ sits on a $b$-triangle, let $w$ and $w'$ be the other extremities of the $b$-edges ending and starting at $v$, respectively. Then $w\ne w'$ and $\Gamma$ has a (non-isolated) $b$-edge from $w'$ to $w$. The \emph{$(\lambda_{2,1},v,w')$-move} consists in removing from $\Gamma$ vertex $v$ and the adjacent edges (the $a$-loop $\ell$ and two $b$-edges). The resulting graph $\Delta$ (see Figure~\ref{fig: ta differenciation}) is $\PSL$-cyclically reduced, it has an isolated $b$-edge from $w'$ to $w$ and combinatorial type $\bm\tau  + \bm\lambda_{2,1}$, where $\bm\lambda_{2,1} =  (-1,0,1,-1,0)$.

If instead $v$ does not sit on a $b$-triangle, there exist vertices $w, w'$ such that $v,w,w'$ are pairwise distinct, there is an isolated $a$-edge between $w$ and $w'$, and an isolated $b$-edge between $v$ and $w$ (two directions are possible for that edge). The \emph{$(\lambda_{2,2}, v\to w,w')$-move} (respectively, $(\lambda_{2,2}, v\leftarrow w,w')$, depending on the orientation of the $b$-edge adjacent to $v$) consists in deleting from $\Gamma$ the vertices $v$ and $w$ and the edges adjacent to them, and adding an $a$-loop $\ell'$ at $w'$. The resulting graph $\Delta$ (see Figure~\ref{fig: ta differenciation}) is $\PSL$-cyclically reduced and has combinatorial type $\bm\tau  + \bm\lambda_{2,2}$, where $\bm\lambda_{2,2} =  (-2,-1,-1,0,0)$.

\begin{lemma}\label{lm: claim 2}
Suppose that $n \ge 3$ and $\ella > 0$.
\begin{enumerate}
\item
The $\lambda_{2,1}$-moves establish a bijection from the set of pairs $(\Gamma,\ell)$ with $\Gamma$ a labeled $\PSL$-cyclically reduced graph of combinatorial type $\bm\tau$ and $\ell$ an $a$-loop adjacent to a $b$-triangle in $\Gamma$, to the set of triples $(\Delta,e,v)$ formed by a labeled $\PSL$-cyclically reduced graph $\Delta$ with combinatorial type $\bm\tau + \bm\lambda_{2,1}$, an isolated $b$-edge $e$ in $\Delta$ and an integer $v\in [n]$.

\item
Similarly, the $\lambda_{2,2}$-moves establish a bijection from the set of pairs $(\Gamma,\ell)$ with $\Gamma$ a labeled $\PSL$-cyclically reduced graph of combinatorial type $\bm\tau$ and $\ell$ an $a$-loop adjacent to an isolated $b$-edge in $\Gamma$, to the set of 4-tuples $(\Delta,\ell',v,w,\epsilon)$ formed by a labeled $\PSL$-cyclically reduced graph $\Delta$ with combinatorial type $\bm\tau + \bm\lambda_{2,2}$, an $a$-loop $\ell'$ in $\Delta$, distinct integers $v,w\in [n]$ and some $\epsilon \in \{-1, +1\}$.
\end{enumerate}
\end{lemma}

\proof
Given a pair $(\Gamma,\ell)$ with $\Gamma$ a labeled $\PSL$-cyclically reduced graph of combinatorial type $\bm\tau$ and $\ell$ an $a$-loop in $\Gamma$ adjacent to a $b$-triangle, we associate to it the triple $(\Delta,e,v)$ --- and we write $(\Gamma,\ell) \longmapsto (\Delta,e,v)$ --- defined as follows: $v$ is the vertex carrying the loop $\ell$ in $\Gamma$; since $v$ is adjacent to a $b$-triangle and we let $e$ be the $b$-edge in that triangle not adjacent to $v$ (going from $w'$ to $w$); finally, we let $\Delta'$ be the weakly labeled graph obtained from $\Gamma$ by a $(\lambda_{2,1},v,w')$-move, $\Delta = \relab(\Delta')$ and $e$ be the isolated $b$-edge in $\Delta$ starting at the vertex labeled $w'$ in $\Delta'$.

Conversely, let $\Delta$ be a labeled $\PSL$-cyclically reduced graph with combinatorial type $\bm\tau + \bm\lambda_{2,1}$, let $e$ be an isolated $b$-edge in $\Delta$ and let $v \in [n]$. Let $\Delta'$ be the weakly labeled graph obtained from $\Delta$ by ``making space for $v$'', that is, by incrementing the labels of all the vertices greater than or equal to $v$. Now let $\Gamma$ be the graph obtained from $\Delta'$ by adding a new vertex $v$, completing $e$ to a $b$-triangle through vertex $v$: it is directly verified that $\Gamma$ is properly labeled, of combinatorial type $\bm\tau$, and that $(\Gamma,\ell) \longmapsto (\Delta,e,v)$. This completes the proof of the first statement.

The second statement is proved in a similar fashion. Given a pair $(\Gamma,\ell)$ with $\Gamma$ a labeled $\PSL$-cyclically reduced graph of combinatorial type $\bm\tau$ and $\ell$ an $a$-loop adjacent to an isolated $b$-edge in $\Gamma$, we associate with it a 4-tuple $(\Delta,\ell',v,w,\epsilon)$ as in the statement, where $v$ is the vertex carrying $\ell$, $w$ is the other extremity of the adjacent isolated $b$-edge and $\epsilon$ records whether a $(\lambda_{2,2}, v\to w,w')$-move or a $(\lambda_{2,2}, v\leftarrow w,w')$ can be performed. The converse mapping, reconstructing $(\Gamma,\ell)$ from $(\Delta,\ell',v,w,\epsilon)$ follows the same steps as for $\lambda_{2,1}$- or $\lambda_3$-moves.
\eop

\begin{figure}[htbp]
\centering
\begin{picture}(130,22)(-2,-22)
\gasset{Nw=4,Nh=4}
\drawoval(37,-12,40,20,12)
\node(n0)(21.0,-07.0){$v'$}
\node(n00)(21.0,-17.0){$w'$}
\node(n1)(4.0,-07.0){$v$}
\node(n11)(4.0,-17.0){$w$}
\drawedge[ELside=r](n1,n11){$b$}
\drawedge[AHnb=0](n1,n0){$a$}
\drawedge[AHnb=0](n00,n11){$a$}

\drawoval(110,-12,40,20,12)
\node(n2)(94.0,-07.0){$v'$}
\node(n22)(94.0,-17.0){$w'$}
\drawedge[AHnb=0](n2,n22){$a$}
\end{picture}
\caption{\small $(\kappa_3,v\to w,v',w')$-move}\label{fig: xb differenciation}
\end{figure}
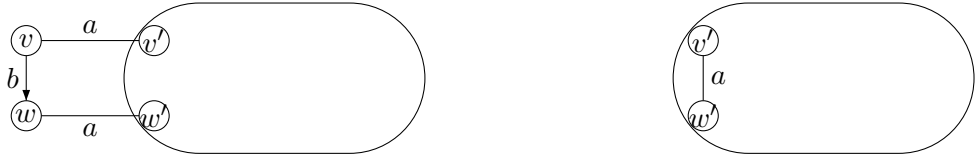

\paragraph{$\kappa_3$-moves}
Let $\Gamma$ be a weakly labeled $\PSL$-cyclically reduced graph of size at least 4, and let $v,w,v',w'$ be pairwise distinct vertices such that there is an isolated $b$-edge from $v$ to $w$, and isolated $a$-edges connecting $v$ and $v'$ on the one hand, and $w$ and $w'$ on the other. The \emph{$\kappa_3$-move} $(\kappa_3, v\to w, v', w')$ consists in deleting vertices $v$ and $w$ and the adjacent edges, and adding a new isolated $a$-edge between $v'$ and $w'$. The resulting graph $\Delta$ (see Figure~\ref{fig: xb differenciation}) is $\PSL$-cyclically reduced and has combinatorial type $\bm\tau  + \bm\kappa_3$, where $\bm\kappa_3 =  (-2,-1,-1,0,0)$.

Similarly to the other moves, we record the following lemma.

\begin{lemma}\label{lm: claim 3}
Suppose that $n\ge 4$, $\ella = 0$ and $\kb > 0$. The $\kappa_3$-moves establish a bijection from the set of pairs $(\Gamma,e)$ with $\Gamma$ a labeled $\PSL$-cyclically reduced graph of combinatorial type $\bm\tau$ and $e$ an isolated $b$-edge, to the set of triples $(\Delta,e',v,w,\epsilon)$ formed by a labeled $\PSL$-cyclically reduced graph $\Delta$ with combinatorial type $\bm\tau + \bm\kappa_{3}$, an isolated $a$-edge $e'$ in $\Delta$, distinct integers $v,w\in [n]$ and some $\epsilon \in \{-1, +1\}$.
\end{lemma}

\proof
The assumption that $\ella = 0$ guarantees that every isolated $b$-edge is adjacent to two isolated $a$-edges, and the fact that $n \ge 4$ guarantees that these $a$-edges are distinct.

Now let $(\Gamma,e)$ be a pair formed by a labeled $\PSL$-cyclically reduced graph $\Gamma$ of combinatorial type $\bm\tau$, and an isolated $b$-edge $e$ in $\Gamma$, say, from vertex $v$ to vertex $w$. Let $v'$ and $w'$ be the other extremities of the isolated $a$-edges adjacent to $v$ and $w$, respectively. We associate with it the tuple $(\Delta,e',v,w,\epsilon)$ --- and we write $(\Gamma,e) \longmapsto (\Delta,e',v,w,\epsilon)$ --- where $\Delta'$ is the weakly labeled graph obtained from $\Gamma$ by a $(\kappa_3,v\to w,v',w')$-move, $\Delta = \relab(\Delta')$, $e'$ is the isolated $a$-edge in $\Delta$ adjacent to the vertex labeled $v'$ in $\Delta'$, $\epsilon = 1$ if $v' < w'$ and $\epsilon = -1$ if $v' > w'$.

Conversely, let $\Delta$ be a labeled $\PSL$-cyclically reduced graph with combinatorial type $\bm\tau + \bm\kappa_3$, let $e'$ be an isolated $a$-edge in $\Delta$, connecting vertices $v'$ and $w'$, and let $v,w$ be distinct integers in $[n]$. Let $\Delta'$ be the weakly labeled graph obtained from $\Delta$ by ``making space for $v,w$'', that is, by incrementing the labels of all the vertices greater than or equal to $\max(v,w)-1$ by 2 units, and those in $[\min(v,w),\max(v,w)-2]$ by 1 unit. Let then $\Gamma$ be the graph obtained from $\Delta'$ by deleting the $a$-edge $e'$; adding new vertices $v,w$ and a $b$-edge from $v$ to $w$; and adding $a$-edges between $v$ and $\min(v',w')$ and between $w$ and $\max(v',w')$ if $\epsilon = 1$ --- between $v$ and $\max(v',w')$ and between $w$ and $\min(v',w')$ if $\epsilon = -1$. It is directly verified that $\Gamma$ is properly labeled, of combinatorial type $\bm\tau$, and that $(\Gamma,e) \longmapsto (\Delta,e',v,w,\epsilon)$.
\eop

\paragraph{Exceptional moves}
Finally, we introduce three so-called \emph{exceptional} moves. The first can be applied only to a weakly labeled version of the 1-vertex graph $\Delta_1$ that does not use label 1, turning it into $\Delta_1$ properly labeled.

The second can be applied only to a weakly labeled version of $\Delta_3$, turning it into $\Delta_1$ (with its only vertex labeled 1). This move can be seen as a degenerate version of a $\lambda_{2,2}$-move. Note that it modifies the combinatorial type by the addition of $\bm{exc} = (-1,0,-1,-1,1)$, the difference between the combinatorial types of $\Delta_1$ and $\Delta_3$.

The last exceptional move can be applied to any weakly labeled version of $\Delta_2$ different from the so-called preferred labeling (see Example~\ref{ex: labelings of small graphs}), turning it to that preferred labeling.

\subsection{Silhouette graphs}

We can see the moves described in Section~\ref{sec: moves} as a rewriting system on weakly labeled $\PSL$-cyclically reduced graphs, which we show to be confluent (Section~\ref{sec: silhouette} below). The following definition will be convenient: we say that a $\PSL$-cyclically reduced graph $\Gamma$ (weakly labeled or not) is a \emph{silhouette graph} if it is equal to $\Delta_1$ or $\Delta_2$, or if it has combinatorial type $(n, n/2, 0, 0, 0)$ (where $n$ is a positive multiple of 6, see Proposition~\ref{prop: charact free and findex}). Observe that no move is defined on a silhouette graph of size at least 3, whichever way it is weakly labeled.

Silhouette graphs play a foundational role in the recursive process for the random generation of subgroups of $\PSL$ described in Section~\ref{sec: random generation}. They also play a central role in \cite{2023:BassinoNicaudWeil-2}.

\subsection{Silhouetting a labeled $\PSL$-cyclically reduced graph}\label{sec: silhouette}

In general, several moves can be applied to a weakly labeled $\PSL$-cyclically reduced graph $\Gamma$. Our next proposition states that, however, the end result of a maximal sequence of moves is independent of the choice of that maximal sequence.

\begin{proposition}\label{prop: uniqueness silhouette}
Let $\Gamma$ be a weakly labeled $\PSL$-cyclically reduced graph $\Gamma$. If $\Delta$ and $\Delta'$ are weakly labeled graphs obtained from $\Gamma$, respectively, after maximal sequences of $\lambda_3$-, $\lambda_{2,1}$-, $\lambda_{2,2}$- and $\kappa_{3}$- and exceptional moves, then $\Delta = \Delta'$.
\end{proposition}

\proof
We proceed by induction on the number of vertices of $\Gamma$. The result is immediate if no move is possible on $\Gamma$. If $\Gamma$ has 1 or 2 vertices, either no move is possible, or only one exceptional move is possible, and the result is again immediate.

Suppose now that $\Gamma$ has $n \ge 3$ vertices and that the sequences of moves leading from $\Gamma$ to $\Delta$ and $\Delta'$ start with the same move (the same type of move, with the same parameters), taking $\Gamma$ to $\Gamma'$. Since $\Gamma'$ is a weakly labeled graph with less than $n$ vertices, the announced result holds by induction.

Finally, suppose that the first moves from $\Gamma$ to $\Delta$ and $\Gamma$ to $\Delta'$, say $m$ and $m'$, are distinct. Note that since $\Gamma$ has size at least 3, neither $m$ nor $m'$ is an exceptional move. Let $\Gamma_1$ (resp. $\Gamma'_1$) be the weakly labeled graph obtained from $\Gamma$ after the move $m$ (resp. $m'$), so that there exists a maximal sequence of moves from $\Gamma_1$ to $\Delta$ (resp. $\Gamma'_1$ to $\Delta'$). Let us consider the possible values of $m$ and $m'$.

If $m = (\lambda_3,v,w)$ and $m' = (\lambda_3,v',w')$, there are two possibilities. If $\Gamma$ is a weakly labeled version of $\Delta_4$ (so that $v = w'$ and $v' = w$), then $\Gamma_1$ and $\Gamma'_1$ are weakly labeled versions of $\Delta_1$, so $\Delta = \Gamma_1$, $\Delta' = \Gamma'_1$ and their normalizations are equal. If instead $\Gamma$ is not a weakly labeled version of $\Delta_4$, then $v, v', w, w'$ are pairwise distinct, and the moves $m$ and $m'$ \emph{commute} in the following sense: an $m'$-move is possible on $\Gamma_1$, leading to a weakly labeled graph $\Gamma_2$; an $m$-move is possible on $\Gamma'_1$, leading to a weakly labeled graph $\Gamma'_2$, and $\Gamma_2 = \Gamma'_2$. Let $\Delta''$ be the graph obtained from $\Gamma_2$ by a maximal sequence of moves.
Since $\Gamma_1$ and $\Gamma'_1$ are weakly labeled graphs with $n-1$ vertices, the induction hypothesis shows that $\Delta = \Delta''$ and $\Delta' = \Delta''$, so that $\Delta = \Delta'$.

We now verify that, similarly, other combinations of first moves $m$ and $m'$ also commute, leading to the same conclusion that the statement in the proposition holds, except in a few degenerate cases that yield the same conclusion by other arguments.

It is readily verified that any $\lambda_3$-move commutes with any $\lambda_{2,1}$-move.

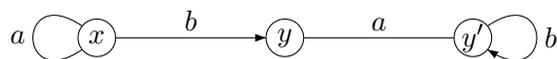
\begin{figure}[htbp]
\centering
\begin{picture}(60,9)(3,-7)
\gasset{Nw=5,Nh=5}
\node(n0)(8.0,-3.0){$x$}
\node(n1)(33.0,-3.0){$y$}
\node(n2)(58.0,-3.0){$y'$}
\drawedge(n0,n1){$b$}
\drawedge[AHnb=0](n1,n2){$a$}
\drawloop[loopangle=180,loopdiam=6,AHnb=0](n0){$a$}
\drawloop[loopangle=0,loopdiam=6](n2){$b$}
\end{picture}
\caption{\small Case where a $\lambda_3$- and a $\lambda_{2,2}$-moves are possible}\label{fig: degenerate 1}
\end{figure}

If $m = (\lambda_3,v,w)$ and $m' = (\lambda_{2,2}, x\rightarrow y,y')$ (or $m' = (\lambda_{2,2}, x\leftarrow y,y')$), we again distinguish two cases. If $y' = v$, then $y = w$ and $\Gamma$ consists of exactly three states $x, y, y'$, with a $b$-edge between $x$ and $y$, an $a$-edge between $y$ and $y'$, an $a$-loop at state $x$ and a $b$-loop at state $y'$, see Figure~\ref{fig: degenerate 1}.
Any maximal sequence of moves from $\Gamma$ leads to $\Delta_1$, where the only vertex is labeled 1, and the announced statement holds. If instead $y' \ne v$, then $\Gamma$ has more than three states and the moves $m$ and $m'$ commute.

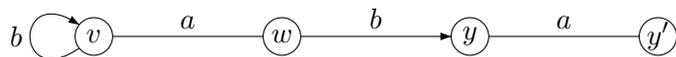
\begin{figure}[htbp]
\centering
\begin{picture}(85,9)(3,-7)
\gasset{Nw=5,Nh=5}
\node(n0)(8.0,-3.0){$v$}
\node(n1)(33.0,-3.0){$w$}
\node(n2)(58.0,-3.0){$y$}
\node(n3)(83.0,-3.0){$y'$}
\drawloop[loopangle=180,loopdiam=6](n0){$b$}
\drawedge[AHnb=0](n0,n1){$a$}
\drawedge(n1,n2){$b$}
\drawedge[AHnb=0](n2,n3){$a$}
\end{picture}
\caption{\small Case where a $\lambda_3$- and a $\kappa_3$-moves are possible}\label{fig: degenerate 2}
\end{figure}

If $m = (\lambda_3,v,w)$ and $m' = (\kappa_3,x\rightarrow y, x', y')$ (or $m' = (\kappa_3,x\leftarrow y, x', y')$), there are again two possibilities. If $v \ne x',y'$, then $m$ and $m'$ modify disjoint parts of $\Gamma$ and they clearly commute. If instead $v = x'$, then $w = x$ --- or if, symmetrically, $v = y'$ and $w = y$ ---, see Figure~\ref{fig: degenerate 2}, then a direct verification shows the following: $m$ can be followed by a $(\lambda_{2,2}, w\rightarrow y, y')$-move (or a $(\lambda_{2,2}, w\leftarrow y, y')$-move, as the case may be), leading to a graph where vertices $v,w,y$ have been deleted and $y'$ carries an $a$-loop; and $m'$ can be followed by a $(\lambda_3,v,y')$-move, leading to that same graph.

The only situation where two $\lambda_2$-moves do not commute is when they are both $\lambda_{2,1}$-moves or both $\lambda_{2,2}$-moves, modifying overlapping parts of $\Gamma$.

\begin{figure}[htbp]
\centering
\begin{picture}(85,20)(3,-19)
\gasset{Nw=5,Nh=5}
\node(n0)(8.0,-3.0){$v$}
\node(n1)(8.0,-19.0){$v'$}
\node(n2)(25.0,-11.0){$w$}
\node(n3)(55.0,-3.0){$v$}
\node(n4)(55.0,-19.0){$v'$}
\node(n5)(72.0,-11.0){$w$}
\node(n6)(92.0,-11.0){$w'$}
\drawloop[loopangle=180,loopdiam=6,AHnb=0](n0){$a$}
\drawloop[loopangle=180,loopdiam=6,AHnb=0](n1){$a$}
\drawloop[loopangle=0,loopdiam=6,AHnb=0](n2){$a$}
\drawloop[loopangle=180,loopdiam=6,AHnb=0](n3){$a$}
\drawloop[loopangle=180,loopdiam=6,AHnb=0](n4){$a$}
\drawedge[ELside=r](n0,n1){$b$}
\drawedge[ELside=r](n1,n2){$b$}
\drawedge[ELside=r](n2,n0){$b$}
\drawedge[ELside=r](n3,n4){$b$}
\drawedge[ELside=r](n4,n5){$b$}
\drawedge[ELside=r](n5,n3){$b$}
\drawedge[AHnb=0](n5,n6){$a$}
\end{picture}
\caption{\small Cases where $\lambda_{2,1}$-moves interfere with each other}\label{fig: 2 lambda21}
\end{figure}
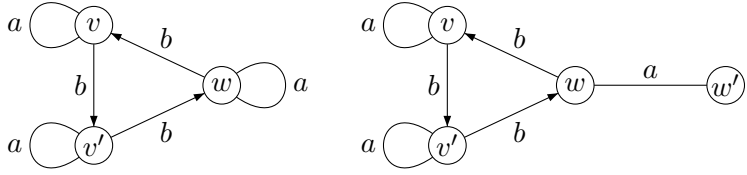

The first case arises if two $a$-loops sit on the same $b$-triangle, so that $m = (\lambda_{2,1},v,w)$ and $m' = (\lambda_{2,1},v',w)$ are possible, see Figure~\ref{fig: 2 lambda21}.
If $w$ also carries an $a$-loop (so that $\Gamma$ has 3 vertices), then $\Delta = \Delta' = \Delta_1$. Otherwise, an isolated $a$-edge connects $w$ and a vertex $w'$, distinct from $v,v',w$, a $(\lambda_{2,2},v'\to w,w)$-move is possible (or some other orientation of a $b$-edge between $v'$ and $w$) after carrying out the $m$-move and, together these two moves amount to deleting vertices $v$, $v'$ and $w$, and adding an $a$-loop at $w'$.

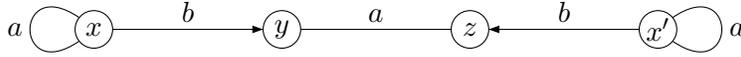
\begin{figure}[htbp]
\centering
\begin{picture}(85,9)(3,-7)
\gasset{Nw=5,Nh=5}
\node(n0)(8.0,-3.0){$x$}
\node(n1)(33.0,-3.0){$y$}
\node(n2)(58.0,-3.0){$z$}
\node(n3)(83.0,-3.0){$x'$}
\drawloop[loopangle=180,loopdiam=6,AHnb=0](n0){$a$}
\drawloop[loopangle=0,loopdiam=6,AHnb=0](n3){$a$}
\drawedge(n0,n1){$b$}
\drawedge[AHnb=0](n1,n2){$a$}
\drawedge[ELside=r](n3,n2){$b$}
\end{picture}
\caption{\small Case where $\lambda_{2,2}$-moves interfere with each other}\label{fig: 2 lambda22}
\end{figure}

In the second case, we have, say, $m = (\lambda_{2,2},x\rightarrow y,z)$ and $m' = (\lambda_{2,2},x'\rightarrow y',z')$ (or any other combination of directions for the $b$-edges between $x$ and $y$ on one hand, and $x'$ and $y'$ on the other) satisfying $y' = z$ and $z' = y$.  Then $\Gamma$ has exactly 4 vertices, see Figure~\ref{fig: 2 lambda22},
and applying either move to $\Gamma$ yields a weakly labeled version of $\Delta_3$, on which one can only apply an exceptional move. It follows that $\Delta = \Delta' = \Delta_1$.

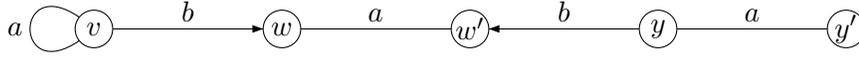
\begin{figure}[htbp]
\centering
\begin{picture}(110,9)(3,-7)
\gasset{Nw=5,Nh=5}
\node(n0)(8.0,-3.0){$v$}
\node(n1)(33.0,-3.0){$w$}
\node(n2)(58.0,-3.0){$w'$}
\node(n3)(83.0,-3.0){$y$}
\node(n4)(108.0,-3.0){$y'$}
\drawloop[loopangle=180,loopdiam=6,AHnb=0](n0){$a$}
\drawedge(n0,n1){$b$}
\drawedge[AHnb=0](n1,n2){$a$}
\drawedge[ELside=r](n3,n2){$b$}
\drawedge[AHnb=0](n3,n4){$a$}
\end{picture}
\caption{\small Case where a $\lambda_{2,2}$- and a $\kappa_3$-moves are possible}\label{fig: degenerate 3}
\end{figure}

It is directly verified that any $\lambda_{2,1}$-move commutes with any $\kappa_3$-move. Consider now the case where $m = (\lambda_{2,2},v\rightarrow w,w')$ and $m'  = (\kappa_3, x\rightarrow y,x',y')$ (or any other combination of directions for the $b$-edges between $v$ and $w$, and between $x$ and $y$). If $w \ne x'$ (and hence $w' \ne x$), then $m$ and $m'$ modify disjoint parts of $\Gamma$ and clearly commute. If instead $w = x'$ and $w' = x$ (see Figure~\ref{fig: degenerate 3}),
a direct verification shows that after applying either $m$ or $m'$, a $\lambda_{2,2}$-move can be applied (a $(\lambda_{2,2},x\rightarrow y,y')$-move can be applied to $\Gamma_1$ and a $(\lambda_{2,2},v\rightarrow w,y')$-move to $\Gamma'_1$), leading to the same weakly labeled graph with an $a$-loop at vertex $y'$: again $m$ and $m'$ commute.

\begin{figure}[htbp]
\centering
\begin{picture}(135,32)(3,-30)
\gasset{Nw=5,Nh=5}
\node(n0)(8.0,-3.0){$v'$}
\node(n1)(33.0,-3.0){$v$}
\node(n2)(58.0,-3.0){$w$}
\node(n3)(83.0,-3.0){$w'$}
\node(n4)(108.0,-3.0){$y$}
\node(n5)(133.0,-3.0){$y'$}
\drawedge[AHnb=0](n0,n1){$a$}
\drawedge[AHnb=0](n2,n3){$a$}
\drawedge[ELside=r](n2,n1){$b$}
\drawedge(n3,n4){$b$}
\drawedge[AHnb=0](n4,n5){$a$}
\node(n10)(58.0,-14.0){$v'$}
\node(n11)(83.0,-14.0){$v$}
\node(n12)(83.0,-30.0){$w$}
\node(n13)(58.0,-30.0){$w'$}
\drawedge[AHnb=0](n10,n11){$a$}
\drawedge[AHnb=0,ELside=r](n12,n13){$a$}
\drawedge[ELside=r](n12,n11){$b$}
\drawedge(n13,n10){$b$}
\end{picture}
\caption{\small Cases where two $\kappa_3$-moves are possible}\label{fig: degenerate 4}
\end{figure}
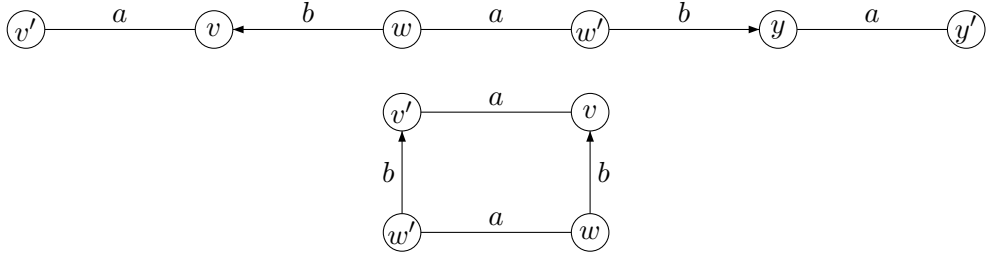

Similarly, suppose that $m = (\kappa_3,v\rightarrow w,v',w')$ and $m' = (\kappa_3,x\rightarrow y,x',y')$ (or any other combination of directions for the $b$-edges). We distinguish three cases. If $x,x',y,y' \notin \{v,v',w,w'\}$, then $m$ and $m'$ clearly commute. If $(w,w') = (x',x)$ and $y \ne v'$ (see Figure~\ref{fig: degenerate 4}), then a $(\kappa_3,x\rightarrow y,v',y')$-move can be applied to $\Gamma_1$ and a $(\kappa_3,v\rightarrow w,v',y')$-move to $\Gamma'_1$, leading to the same weakly labeled graph (with an $a$-edge between $v'$ and $y'$, and no vertices labeled $x, y, v, w$).
If now $(w,w') = (x',x)$ and $y = v'$ (so that $y' = v$), then $\Gamma$ has exactly four vertices, both $\Gamma_1$ and $\Gamma'_1$ are weak labelings of $\Delta_2$, on which only an exceptional move is defined, so that $\Delta = \Delta'$.

This concludes the proof of the proposition.
\eop

If $\Gamma$ is a labeled $\PSL$-cyclically reduced graph, we define the \emph{silhouette} $\core(\Gamma)$ of $\Gamma$ to be the labeled graph resulting from the application of a maximal sequence of moves, followed by a normalization: Proposition~\ref{prop: uniqueness silhouette} guarantees that $\core(\Gamma)$ is well defined.

\begin{example}\label{ex: silhouetting}
Consider the three (labeled) Stallings graphs in Figure~\ref{fig: stallings graphs}. The first is equal to its own silhouette, and is also equal to the silhouette of the third. The silhouette of the second graph is $\Delta_2$ (with its preferred labeling).
\end{example}

The silhouetting operation preserves some important algebraic information about a subgroup, namely the free rank component of its isomorphism type.

\begin{proposition}\label{prop: rank and small silhouette}
Let $H$ be a $\PSL$-cyclically reduced subgroup of $\PSL$, with Stallings graph $\Gamma$ and isomorphism type $(\ella, \ellb, r)$. If $\core(\Gamma)$ has isomorphism type $(\ella',\ellb',r')$, then $r = r'$. In particular $\core(\Gamma) = \Delta_1$ (respectively, $\Delta_2$) if and only if $r = 0$ (respectively, $r = 1$).
\end{proposition}

\proof
Let $\bm\tau$ be the combinatorial type of $\Gamma$. Proposition~\ref{prop: combinatorial vs. isomorphism} shows that the free rank $r$ in the isomorphism type of $H$ is a function of $\bm\tau$; more precisely, if $\bm\tau = (n,\ka,\kb,\ella,\ellb)$, then $6(r-1) = n-2\kb-3\ella-4\ellb = \phi(\bm\tau)$, and we observe that $\phi$ is a linear map.

By construction, $\core(\Gamma)$ is obtained from $\Gamma$ by a succession of $\lambda_3$-, $\lambda_{2,1}$-, $\lambda_{2,2}$-, $\kappa_3$-moves and maybe one  exceptional move (followed by normalization). Each of these moves modifies the combinatorial type by adding to it the vector $\bm\lambda_3$, $\bm\lambda_{2,1}$, $\bm\lambda_{2,2}$, $\bm\kappa_3$ or $\bm{exc}$. Every one of these vectors lies in the kernel of $\phi$, so the free rank component of the isomorphism types of $\Gamma$ and $\core(\Gamma)$ coincide.

It is immediate that this free rank component is 0 for $\Delta_1$, 1 for $\Delta_2$ and $1 + n/6 \ge 2$ for each silhouette graph of size $n > 2$. The proposition follows immediately.
\eop

\section{Counting subgroups by isomorphism and by combinatorial type}\label{sec: counting}

Our aim in this section is to count subgroups of a given size, under some additional constraint: with a fixed isomorphism type or with a fixed combinatorial type. Since each subgroup is uniquely represented by its Stallings graph, \emph{i.e.}, by a $\PSL$-reduced graph, this is equivalent to counting these graphs (up to isomorphism).

As noted in Section~\ref{sec: labeled graphs}, an $n$-vertex $\PSL$-reduced graph admits exactly $n!$ distinct labelings. As a result, our strategy to count $n$-vertex $\PSL$-reduced graphs will be to count labeled $n$-vertex $\PSL$-reduced graphs, and then divide that count by $n!$. The same applies for the counting of $n$-vertex $\PSL$-reduced graphs of a particular combinatorial type, or for $n$-vertex rooted $\PSL$-cyclically reduced graphs. Note that there is no such easy correlation between the number of labeled and unlabeled (non-rooted) cyclically reduced graphs, as counting is perturbed by the existence of symmetries (automorphisms).

Thus our task reduces to counting labeled $\PSL$-reduced graphs. It further reduces to counting labeled $\PSL$-cyclically reduced graphs, as we explain below.

\subsection{Reduction to the count of labeled $\PSL$-cyclically reduced graphs}\label{sec: reduction to cyclically reduced}

If $\bm\tau = (n,\ka,\kb,\ella,\ellb)$ is a tuple of integers, we let $H(\bm\tau)$ (respectively, $L(\bm\tau)$, $s(\bm\tau)$) be the number of subgroups (respectively, labeled $\PSL$-reduced graphs, labeled $\PSL$-cyclically reduced graphs) of combinatorial type $\bm\tau$.

\begin{example}\label{ex: S coefficients small n}
In view of Example~\ref{ex: small n}, the non-zero values of $H$, $L$ and $s$ for tuples $(n,\ka,\kb,\ella,\ellb)$ where $n = 1,2$ are as follows:

$\bullet$\enspace $H(\bm\tau)  = L(\bm\tau) = 1$ for $\bm\tau = (1,0,0,1,1), (1,0,0,1,0), (1,0,0,0,1)$ and $s(\bm\tau) = 1$ for $\bm\tau = (1,0,0,1,1)$;

$\bullet$\enspace $L(\bm\tau) = 4$ and $H(\bm\tau) = s(\bm\tau) = 2$ for $\bm\tau = (2,1,1,0,0), (2,0,1,2,0)$; 

$\bullet$\enspace $L(\bm\tau) = 2$ and $H(\bm\tau) = s(\bm\tau) = 1$ for $\bm\tau = (2,1,0,0,2)$; 

$\bullet$\enspace $L(\bm\tau) = 2$ and $H(\bm\tau) = 1$ for $\bm\tau = (2,0,1,1,0), (2,1,0,0,1)$; 
\end{example}

We first establish the connection between the parameters $H(\bm\tau)$, $L(\bm\tau)$ and $s(\bm\tau)$.

\begin{proposition}\label{prop: counting subgroups by cb type}
Let $\bm\tau = (n,\ka,\kb,\ella,\ellb)$ be a combinatorial type with $n\ge 2$. The numbers $H(\bm\tau)$, $L(\bm\tau)$ and $s(\bm\tau)$, respectively of subgroups, labeled $\PSL$-reduced graphs and labeled $\PSL$-cyclically reduced graphs of combinatorial type $\bm\tau$ are related as follows.
$$L(\bm\tau) = n\cdot s(n,\ka,\kb,\ella,\ellb) + (\ella+1)\cdot s(n,\ka,\kb,\ella+1,\ellb) + (\ellb+1)\cdot s(n,\ka,\kb,\ella,\ellb+1) $$
$$H(\bm\tau) = \frac1{n!}\,L(\bm\tau).$$
\end{proposition}

\proof 
Let $(\Gamma,v)$ be a $\PSL$-reduced graph with $n\ge 2$ vertices, such that $\Gamma$ is not $\PSL$-cyclically reduced. Then $v$ is adjacent to an $a$-edge but no $b$-edge, or the opposite. Adding a $b$-loop at $v$ in the first case, an $a$-loop in the second case, yields a rooted $\PSL$-cyclically reduced graph $(\Gamma',v)$. Conversely, if $\Gamma'$ is $\PSL$-cyclically reduced, we get $\PSL$-reduced graphs either by rooting $\Gamma'$ at any one of its vertices, or by rooting $\Gamma'$ at a vertex that carries a loop and deleting that loop. The first equality follows directly.

The second equality follows from the first since a size $n$ $\PSL$-reduced graph has $n!$ distinct labelings (see Section~\ref{sec: labeled graphs}).
\eop

Based on Proposition~\ref{prop: combinatorial vs. isomorphism}, which relates the isomorphism type and the combinatorial type of a subgroup, we get the following statement.

\begin{proposition}\label{prop: counting subgroups by isom type}
Let $\bm\sigma = (\ella,\ellb,r)$ be an isomorphism type and let $\ka = \frac12(n-\ella)$.

The number of $\PSL$-cyclically reduced subgroups of size $n$ and isomorphism type $\bm\sigma$ is $n\cdot s(n,\ka,\kb,\ella,\ellb)$, where $\kb = \frac12(n - 3\ella - 4\ellb - 6r + 6)$.

The number of non-$\PSL$-cyclically reduced subgroups of size $n$ and isomorphism type $\bm\sigma$, where the base vertex is adjacent to an $a$-edge, is $(\ellb + 1)\cdot s(n,\ka,\kb',\ella,\ellb+1)$, where $\kb' = \frac12(n - 3\ella - 4\ellb - 6r + 2)$.

The number of non-$\PSL$-cyclically reduced subgroups of size $n$ and isomorphism type $\bm\sigma$, where the base vertex is adjacent to a $b$-edge, is $(\ella + 1)\cdot s(n,\ka,\kb'',\ella+1,\ellb)$, where $\kb'' = \frac12(n - 3\ella - 4\ellb - 6r + 4)$.
\end{proposition}

Propositions~\ref{prop: counting subgroups by cb type} and~\ref{prop: counting subgroups by isom type} effectively reduce the counting of subgroups to the counting of labeled $\PSL$-cyclically reduced graphs of a given combinatorial type, which is investigated in Section~\ref{sec: couting cyclically reduced} below.

\subsection{Counting labeled $\PSL$-cyclically reduced graphs}\label{sec: couting cyclically reduced}

Let $\bm\tau = (n,\ka,\kb,\ella,\ellb)$ be a combinatorial type. We give (multi-)recurrence relations to compute $s(\bm\tau)$ when $n > 2$. For $n\le 2$, see Example~\ref{ex: S coefficients small n}.

The bijections established in Lemmas~\ref{lm: claim 1}, \ref{lm: claim 2} and~\ref{lm: claim 3} show the following.

\begin{proposition}\label{prop: claim 4}
Let $\bm\tau = (n,\ka,\kb,\ella,\ellb)$ be a combinatorial type such that $n \ge 2$ and $\ellb > 0$. If $\Delta$ has combinatorial type $\bm\tau + \bm\lambda_3$, then the set of labeled $\PSL$-cyclically reduced graphs $\Gamma$ of combinatorial type $\bm\tau$, such that a $\lambda_3$-move takes $\Gamma$ to $\Delta$, has $\frac{n\cdot (\ella + 1)}\ellb$ elements. More generally,
\begin{align}
s(\bm\tau) &= \frac{n\cdot (\ella + 1)}\ellb\ s(\bm\tau+\bm\lambda_3),\textrm{ that is:}\nonumber\\
s(n,\ka,\kb,\ella,\ellb) &= \frac{n\cdot (\ella + 1)}\ellb\ s(n-1,\ka-1,\kb,\ella+1,\ellb-1).\label{eq: fromStb}
\end{align}
\end{proposition}

\begin{proposition}\label{prop: claim 5}
Let $\bm\tau = (n,\ka,\kb,\ella,\ellb)$ be a combinatorial type such that $n \ge 3$ and $\ella > 0$. If $\Delta$ has combinatorial type $\bm\tau + \bm\lambda_{2,1}$ (resp. $\bm\tau + \bm\lambda_{2,2}$), then the set of labeled $\PSL$-cyclically reduced graphs $\Gamma$ of combinatorial type $\bm\tau$, such that a $\lambda_{2,1}$-move (resp. a $\lambda_{2,2}$-move) takes $\Gamma$ to $\Delta$, has $\frac{n\cdot (\kb + 1)}\ella$ (resp. $2n\cdot (n-1)$) elements. More generally,
\begin{align}
s(\bm\tau) \enspace=\enspace &\frac{n\cdot (\kb + 1)}{\ella}\ s(\bm\tau+\bm\lambda_{2,1})\nonumber \\
& + 2n\cdot (n-1)\ s(\bm\tau+\bm\lambda_{2,2}),\textrm{ that is:}\nonumber\\
s(n,\ka,\kb,\ella,\ellb) \enspace=\enspace &\frac{n\cdot (\kb + 1)}{\ella}\ s(n-1,\ka,\kb+1,\ella-1,\ellb)\nonumber \\
& + 2n\cdot (n-1)\ s(n-2,\ka-1,\kb-1,\ella,\ellb).\label{eq: fromSta}
\end{align}
\end{proposition}

\begin{proposition}\label{prop: claim 6}
Let $\bm\tau = (n,\ka,\kb,\ella,\ellb)$ be a combinatorial type such that $n \ge 4$, $\ella = 0$ and $\kb > 0$. If $\Delta$ has combinatorial type $\bm\tau + \bm\kappa_3$, then the set of labeled $\PSL$-cyclically reduced graphs $\Gamma$ of combinatorial type $\bm\tau$, such that a $\kappa_3$-move takes $\Gamma$ to $\Delta$, has $\frac{2n\cdot (n-1)(\ka - 1)}\kb$ elements. More generally, if $\ella = 0$, we have
\begin{align}
s(\bm\tau) \enspace=\enspace & 2\ \frac{n\cdot (n-1)(\ka - 1)}{\kb}\ s(\bm\tau + \bm\kappa_3),\textrm{ that is:}\nonumber\\
s(n,\ka,\kb,0,\ellb) \enspace=\enspace & 2\ \frac{n\cdot (n-1)(\ka - 1)}{\kb}\ s(n-2,\ka-1,\kb-1,0,\ellb).\label{eq: fromSy}
\end{align}
\end{proposition}

We can use Equations~(\ref{eq: fromStb}), (\ref{eq: fromSta}) and~(\ref{eq: fromSy}) to compute the coefficient $s(n,\ka,\kb,\ella,\ellb)$, where $n \ge 3$: if one of $\kb$, $\ella$ or $\ellb$ is greater than zero, we can apply at least one of these equations, thus reducing the first argument of the coefficients to compute by 1 or 2.

More precisely, one may first iterate the use of Equation~(\ref{eq: fromStb}) until $n \le 2$ or $\ellb = 0$. One can then use repeatedly Equation~(\ref{eq: fromSta}), thus reducing the computation of $s(\bm\tau)$ to the computation of a number of smaller coefficients, until $n\le 2$ or $\ella = 0$ (note that Equation~(\ref{eq: fromSta}) never increases $\ellb$). Finally, if $n \ge 3$ and $\ella = \ellb = 0$, then in fact $n\ge 4$ and one can use repeatedly Equation~(\ref{eq: fromSy}) until $n \le 2$ or $\kb = \ella = \ellb = 0$. The computation of the coefficients when $n \le 2$ was done in Example~\ref{ex: S coefficients small n}. As for the coefficients of the form $s(n,\ka,0,0,0)$ ($n > 2$), we note that they count the size $n$ labeled silhouette graphs.

The latter numbers were computed in \cite[Appendices A.3 and A.4]{2021:BassinoNicaudWeil} (see also the computation by Stothers \cite{1978:Stothers-MathsComp} of the number of finite index, free subgroups of $\PSL$, that is, of subgroups whose Stallings graph is a silhouette graph of size at least 3). 

\begin{proposition}\label{prop: counting silhouettes}
Let $t_2$ (respectively, $t_3$) be given, for $n\ge 1$, by\footnote{What is written $t_2(2n)$ (respectively, $t_3(3n)$, $s(6n,3n,0,0)$) here, is written $t_2^{(0)}(2n)$ (respectively, $t_3^{\textsf{fr-fi}}(3n)$, $g_{\textsf{pr}}^{\textsf{fr-fi}}(6n)$) in \cite{2021:BassinoNicaudWeil}.}
\begin{align*}
t_2(2n) = \frac{(2n)!}{2^n\, n!} = \prod_{1\le i\le n} (2i-1), \quad
\textrm{and }& t_3(3n) = \frac{(3n)!}{3^n\, n!} = \prod_{1\le i\le n} (3i-1)(3i-2). 
\end{align*}
Then the number $s(6n,3n,0,0,0)$ of size $6n$ labeled silhouette graphs ($n \ge 1$) satisfies the following recurrence relation:
\begin{align*}
s(6n,3n,0,0,0) &= t_2(6n)\,t_3(6n) - \sum_{m=1}^{n-1} t_2(6m)\,t_3(6m)\,s\big(6(n-m),3(n-m),0,0,0\big).
\end{align*}
\end{proposition}

\section{Random generation of subgroups of $\PSL$}\label{sec: random generation}

Our objective in this section is to produce an algorithm which generates uniformly at random subgroups of $\PSL$ with a given size and isomorphism type. 

As we saw in Example~\ref{ex: S coefficients small n}, there are exactly four size 1 subgroups, with pairwise distinct combinatorial and isomorphism type: the trivial subgroup, the subgroups generated by $a$ and $b$, respectively, and $\PSL$ itself. We now concentrate on generating subgroups of size at least $2$, and we assume that the parameters $L(\bm\tau)$ and $s(\bm\tau)$ have been pre-computed for all types of sufficient size.

Like in Section~\ref{sec: counting}, generating uniformly at random a subgroup of a given combinatorial or isomorphism type reduces to randomly generating a labeled $\PSL$-reduced graph of a given combinatorial type and, before that, to randomly generating a labeled $\PSL$-cyclically reduced graph of a given type. Indeed, the label-forgetting map, from the set of labeled $\PSL$-reduced graphs of combinatorial type $\bm\tau = (n,\ka,\kb,\ella,\ellb)$ to the set of $\PSL$-reduced graphs of type $\bm\tau$, is such that the inverse image of each $\PSL$-reduced graph of type $\bm\tau$ contains exactly $n!$ elements (see the discussion at the end of Section~\ref{sec: labeled graphs}).

As we saw in Proposition~\ref{prop: combinatorial vs. isomorphism}, the isomorphism class of a $\PSL$-reduced graph is determined by its combinatorial type, and a given size and isomorphism type arises for a finite number of combinatorial types only. As a result, we only need to randomly generate a $\PSL$-reduced graph with a given combinatorial type, and this starts with randomly generating a labeled $\PSL$-cyclically reduced graphs of a given combinatorial type.

We first deal with the particular case of labeled silhouette graphs, then proceed to the general case of labeled $\PSL$-cyclically reduced graphs and, finally, to labeled $\PSL$-reduced graphs.

\subsection{Random labeled silhouette graphs}\label{sec: drawing a silhouette graph}

Let $n$ be a positive multiple of $6$. The procedure $\rsilh(n)$ to generate a size $n$ labeled silhouette graph, summarized below, is well known (see \cite{2021:BassinoNicaudWeil} for instance). If $s$ is a permutation on $n$ elements, we denote by $\shuffle(s)$ the permutation $t\inv s t$ where $t$ is a random permutation on $n$ elements.

\begin{algorithm}[H]
\DontPrintSemicolon

\SetKwInOut{Input}{input}
\SetKwInOut{Output}{output}
\SetKwRepeat{Do}{do}{while}

\Do{
The graph $\Gamma$ determined by $s_2$ and $s_3$ is not connected
}{
$s_2 = \shuffle((1\ 2)\,(3\ 4)\ldots(n-1\ n))$\;
$s_3 = \shuffle((1\ 2\ 3)\,(4\ 5\ 6)\ldots(n-2\ n-1\ n))$\;
}
\caption{$\rsilh(n)$}
\end{algorithm}
\medskip

Note that the random permutations $s_2$ and $s_3$ may well determine a disconnected graph, but the proof of \cite[Proposition 8.18]{2021:BassinoNicaudWeil} shows that this happens with vanishing probability (precisely: with probability $\frac5{36}n\inv + o(n\inv)$). Therefore this algorithm (a rejection algorithm) produces a silhouette graph after $k$ iterations, with $\E(k) \sim 1$.

\subsection{Random $\PSL$-cyclically reduced graphs}\label{sec: drawing a PSL cyclically reduced graph}

We exploit, again, the bijections established in Section~\ref{sec: moves}, which we already used to derive the recurrence relations in Section~\ref{sec: couting cyclically reduced}. This yields the following algorithm, called $\rcrg(\bm\tau)$, to randomly generate a $\PSL$-cyclic\-ally reduced graph of combinatorial type $(\bm\tau)$.

We use the following notation: if $v$ is an integer, $\shift_v$ is the map defined on integers by $\shift_v(x) = x$ if $x< v$ and $\shift_v(x) = x+1$ if $x\ge v$; if $v, w$ are distinct integers, $\shift_{v,w}$
 is the map defined on integers by $\shift_{v,w}(x) = x$ if $x < \min(v, w)$, $\shift_{v,w}(x) = x+1$ if $\min(v,w) \le x < \max(v, w)-1$,  and $\shift_{v,w}(x) = x+2$ if $x \ge \max(v, w)-1$. Note that $\shift_v$ ``pushes'' all integers greater than or equal to $v$ by one unit, so that the range of $\shift_v$ misses $v$; similarly, the range of $\shift_{v,w}$ misses $v$ and $w$.

We extend this notation to any graph $\Delta$ labeled by integers: if $v$ is an integer, the graph $\shift_v(\Delta)$ is a relabeling of the vertices of $\Delta$ using $\shift_v$ on each vertex label; if $v$ and $w$ are two distinct integers, the graph $\shift_{v,w}(\Delta)$ is a relabeling of the vertices of $\Delta$ using $\shift_{v,w}$ on each vertex label.

\begin{algorithm}[H]
\small
\DontPrintSemicolon
\LinesNumbered

\If{$\bm\tau = (1,0,0,1,1)$\label{rcrg:if delta1}}{\Return the unique labeled $\Delta_1$}
\If{$\bm\tau = (2,1,1,0,0)$\label{rcrg:if delta2}}{\Return any one of the two labeled $\Delta_2$}
\If{$\bm\tau = (2,0,1,2,0)$\label{rcrg:if delta3}}{\Return any one of the two labeled $\Delta_3$}
\If{$\bm\tau = (1,0,0,1,1)$\label{rcrg:if delta4}}{\Return the unique labeled $\Delta_4$}
\BlankLine
\tcp{At this stage $n$ is necessarily greater than 2}
\If{$\bm\tau = (n,\ka,\kb,\ella,\ellb)$ and $\ellb>0$\label{rcrg:ellb}}{
	$\Delta = \rcrg(\bm\tau+\bm\lambda_3)$\label{rcrg:delta ellb}\;
	$w=$ uniform random vertex with an $a$-loop $\ell'$ in $\Delta$\; 
	$v=$ uniform random integer in $\{1,\dots, n\}$\; 
	\Return $\Gamma$ constructed from $\Delta$ by relabeling its vertices  using $\shift_v$, removing the $a$-loop $\ell'$ at $\shift_v(w)$, adding a new vertex labeled $v$ and a $b$-loop at $v$, and adding an $a$-edge between $v$ and $\shift_v(w)$\;
}
\If{$\bm\tau = (n,\ka,\kb,\ella,0)$ and $\ella>0$\label{rcrg:ella}}{
	$x=$ uniform integer in $[s(n,\ka,\kb,\ella,0)]$\;\label{rcrg:random x}
	\eIf{$x\leq n\cdot (\kb+1)\cdot s(n-1,\ka,\kb+1,\ella-1,0)$\label{rcrg:inner if}}{\label{rcrg:compare x}
		$\Delta = \rcrg(\bm\tau+\bm\lambda_{2,1})$\;
        $(w\xrightarrow{b} w') =$ uniform random isolated $b$-edge in $\Delta$\;
        $v=$ uniform random  integer in $\{1,\dots, n\}$\;
        \Return $\Gamma$    constructed from $\Delta$ by relabeling its vertices using $\shift_v$, adding a new vertex labeled $v$ and an $a$-loop at $v$, and adding $b$-edges from $\shift_v(w')$ to $v$ and from $v$ to $\shift_v(w)$\;
	}{
		$\Delta = \rcrg(\bm\tau+\bm\lambda_{2,2})$\;
		$w'=$ uniform random vertex with a $a$-loop $\ell'$ in $\Delta$\;
		$(v,w)=$ uniform random pair of distinct integers in $\{1,\dots, n\}$\;
		\Return $\Gamma$ constructed from $\Delta$ by removing the $a$-loop $\ell'$, relabeling the vertices of $\Delta$ using $\shift_{v,w}$, adding new vertices labeled $v$ and $w$, an $a$-edge between $w$ and $\shift_{v,w}(w')$, a $b$-edge between $v$ and $w$ (choosing orientation uniformly at random) and an $a$-loop at $v$\;
	}
}
\If{$\bm\tau = (n,\ka,\kb,0,0)$ and $\kb>0$\label{rcrg:kb}}{
	$\Delta = \rcrg(\bm\tau+\bm\kappa_3)$\;
	$(v'\stackrel{\substack{\\[-2pt] a}}{-} w') =$ uniform random isolated $a$-edge  in $\Delta$\;
	$(v,w)=$ uniform random pair of distinct integers in $\{1,\dots, n\}$\;
	\Return $\Gamma$ constructed from $\Delta$ by removing the the $a$-edge $e'$, relabeling the vertices using $\shift_{v,w}$, adding a $b$-edge between $v$ and $w$ (choosing its orientation uniformly at random) and $a$-edges between $v$ and $\shift_{v,w}(v')$, and between $w$ and $\shift_{v,w}(w')$, respectively\;
}
\tcp{At this stage $\ella=\ellb=\kb=0$}
\Return $\rsilh(n)$\label{rcrg:silhouette}

\caption{$\rcrg(n)$}\label{algo:rcrg}
\end{algorithm}

\begin{theorem}\label{thm: proof of algorithm rcrg}
Algorithm $\rcrg$, described above, produces, on input a combinatorial type $\bm\tau$, a random $\PSL$-cyclically reduced graph of type $\bm\tau$.
\end{theorem}

\proof
Let $\bm\tau = (n,\ka,\kb,\ella,\ellb)$. We proceed by induction on $n$. Note that Algorithm $\rcrg$ is recursive, and that, for any $\bm\tau$, only one of the outer \texttt{if} statements (Lines~\ref{rcrg:if delta1}, \ref{rcrg:if delta2}, \ref{rcrg:if delta3}, \ref{rcrg:if delta4}, \ref{rcrg:ellb}, \ref{rcrg:ella} and \ref{rcrg:kb}) holds. Moreover, the algorithm stops immediately after the \texttt{if} statements of Lines~\ref{rcrg:if delta1}, \ref{rcrg:if delta2}, \ref{rcrg:if delta3}, \ref{rcrg:if delta4}, and strictly decreases the value of $n$ for the other ones. As a result, Algorithm $\rcrg$ stops on any input $\bm\tau$ (which is a proper combinatorial type).

The statement of the theorem holds trivially if $n \le 2$. Let us now assume that $n > 2$. 

If $\ellb > 0$ (that is, if the condition of Line~\ref{rcrg:ellb} holds), Lemma~\ref{lm: claim 1} describes a bijection between the set of pairs $(\Gamma,\ell)$, where $\Gamma$ is a $\PSL$-cyclically reduced graph of combinatorial type $\bm\tau$ and $\ell$ is a $b$-loop in $\Gamma$, and the set of triples $(\Delta,\ell',v)$ where $\Delta$ is a $\PSL$-cyclically reduced graph of type $\bm\tau+\bm\lambda_3$, $\ell'$ is an $a$-loop in $\Delta$ and $v \in [n]$. This bijection, between two finite sets, preserves uniformity. Thus the first steps in this case (selecting uniformly at random $\Delta$, $\ell'$ and $v$) translate into the selection, uniformly at random, of a pair $(\Gamma,\ell)$ where $\Gamma$ has combinatorial type $\bm\tau$ and $\ell$ is one of the $\ellb$ $b$-loops in $\Gamma$ (a number that depends on $\bm\tau$ but not on $\Gamma$). Forgetting the $\ell$-component of this pair yields a randomly chosen $\PSL$-cyclically reduced graph of type $\bm\tau$.

The reasoning is exactly similar if the condition of Line~\ref{rcrg:kb} holds, relying on Lemma~\ref{lm: claim 3}.

For the condition of Line~\ref{rcrg:ella}, we need to handle the two options, corresponding to $\lambda_{2,1}$- and $\lambda_{2,2}$-moves. The set of pairs $(\Gamma,\ell)$, where $\Gamma$ is a $\PSL$-cyclically reduced graph of combinatorial type $\bm\tau$ and $\ell$ is an $a$-loop in $\Gamma$, is partitioned in two subsets $S_1$ and $S_2$: $(\Gamma,\ell)\in S_1$ if $\ell$ is adjacent to a $b$-triangle, and $(\Gamma,\ell)\in S_2$ if $\ell$ is adjacent to an isolated $b$-edge. Lemma~\ref{lm: claim 2} describes the sets $S_1$ and $S_2$ are in bijection with. This determines the cardinalities of $S_1$ and $S_2$, which correspond precisely to the probability tested at Line~\ref{rcrg:inner if}. The reasoning is then identical to Lines~\ref{rcrg:ellb} and~\ref{rcrg:kb}.

Finally, if none of these conditions holds (so that $\bm\tau$ is the combinatorial type of a silhouette graph), the algorithm uses the \texttt{return} command on Line~\ref{rcrg:silhouette} to produce a random silhouette graph, see Section~\ref{sec: drawing a silhouette graph}.
\eop

\subsection{Random subgroups of $\PSL$}\label{sec: drawing a PSL reduced graph}

We show how to randomly generate subgroups of a given combinatorial type, and then of a given size and isomorphism type.

\subsubsection*{Random generation for a given combinatorial type}

Let $\bm\tau = (n,\ka,\kb,\ella,\ellb)$ be a combinatorial type. The formula for the number $L(\bm\tau)$ of labeled $\PSL$-reduced graphs of type $\bm\tau$, in Proposition~\ref{prop: counting subgroups by cb type} above, suggests the following algorithm to draw uniformly at random a labeled $\PSL$-reduced graph of combinatorial type $\bm\tau$.

\begin{itemize}
  \item[(1)] Draw an integer $0 \le p < L(\bm\tau)$ uniformly at random.
  \item[(2)] If $p < n \cdot s(\bm\tau)$
    and $q$ is the quotient of $p$ by $s(\bm\tau)$ (so that $0\le q < n$), draw uniformly at random a labeled $\PSL$-cyclically reduced graph with combinatorial type $\bm\tau$ and root it at vertex $q+1$.
  \item[(3)] If $n\cdot s(\bm\tau) \leq p < n\cdot s(\bm\tau) + (\ella+1)\cdot s(n,\ka,\kb,\ella+1,\ellb)$
    and $q$ is the quotient of $p - n\cdot s(\bm\tau)$ by $s(n,\ka,\kb,\ella+1,\ellb)$ (so that $0\le q \le \ella$), draw uniformly at random a labeled $\PSL$-cyclically reduced graph with combinatorial type $(n,\ka,\kb,\ella+1,\ellb)$ (as in Section~\ref{sec: drawing a PSL cyclically reduced graph}), delete the $(q+1)$st $a$-loop (following the order of vertex labels) and root the graph at the vertex where that loop used to be.
    \item[(4)] If $n\cdot s(\bm\tau) + (\ella+1)\cdot s(n,\ka,\kb,\ella+1,\ellb) \le p$
    and $q$ is the quotient of $p - n\cdot s(\bm\tau) - (\ella+1)\cdot s(n,\ka,\kb,\ella+1,\ellb)$ by $s(n,\ka,\kb,\ella,\ellb+1)$ (so that $0\le q \le \ellb$), draw uniformly at random a labeled $\PSL$-cyclically reduced graph with combinatorial type $(n,\ka,\kb,\ella,\ellb+1)$ (as in Section~\ref{sec: drawing a PSL cyclically reduced graph}), delete the $(q+1)$st $b$-loop (following the order of vertex labels) and root the graph at the vertex where that loop used to be.
\end{itemize}

To draw uniformly at random a subgroup of combinatorial type $\bm\tau$, we first draw a labeled $\PSL$-reduced graph of type $\bm\tau$, and then forget the labeling.

\begin{remark}\label{rk: drawing cyclically reduced by cb type}
  To draw uniformly at random a $\PSL$-cyclically reduced subgroup of combinatorial type $\bm\tau$, the algorithm is modified as follows: in step (1), one draws an integer $p$ between 0 and $n\cdot s(\bm\tau)-1$; one then applies only step (2).
\end{remark}

\subsubsection*{Random generation for a given size and isomorphism type}

Now let $n$ be a positive integer and let $\bm\sigma = (\ella,\ellb,r)$ be an isomorphism type. Let $\ka = \frac12(n-\ella)$, $\kb = \frac12(n - 3\ella - 4\ellb - 6r + 6)$, $\kb' = \frac12(n - 3\ella - 4\ellb - 6r + 2)$ and $\kb'' = \frac12(n - 3\ella - 4\ellb - 6r + 4)$.

Proposition~\ref{prop: counting subgroups by isom type} suggests the following algorithm to draw uniformly at random a subgroup of size $n$ and isomorphism type $\bm\sigma$.

\begin{itemize}
\item[(1)] Draw uniformly at random an integer $p$ between 0 and
$$n\cdot s(n,\ka,\kb,\ella,\ellb) + (\ellb + 1)\cdot s(n,\ka,\kb',\ella,\ellb+1) + (\ella + 1)\cdot s(n,\ka,\kb'',\ella+1,\ellb)-1.$$
\item[(2)] If $p < n\cdot s(n,\ka,\kb,\ella,\ellb)$,
  draw uniformly at random a labeled rooted $\PSL$-cyclically reduced graph with combinatorial type $(n,\ka,\kb,\ella,\ellb)$.
\item[(3)] If $n\cdot s(\bm\tau) \le  p < n\cdot s(\bm\tau) + (\ellb+1)\cdot s(n,\ka,\kb',\ella,\ellb+1)$
and $q$ is the quotient of $p - n\cdot s(\bm\tau)$ by $s(n,\ka,\kb',\ella,\ellb+1)$ (so that $0\le q \le \ellb$), draw uniformly at random a labeled $\PSL$-cyclically reduced graph with combinatorial type $(n,\ka,\kb',\ella,\ellb+1)$, delete the $(q+1)$st $b$-loop (following the order of vertex labels) and root the graph at the vertex where that loop used to be.
\item[(4)] If $n\cdot s(\bm\tau) + (\ellb+1)\cdot s(n,\ka,\kb',\ella,\ellb+1) \le p$
and $q$ is the quotient of $p - n\cdot s(\bm\tau) - (\ellb+1)\cdot s(n,\ka,\kb',\ella,\ellb+1)$ by $s(n,\ka,\kb'',\ella+1,\ellb)$ (so that $0\le q \le \ella$), draw uniformly at random a labeled $\PSL$-cyclically reduced graph with combinatorial type $(n,\ka,\kb,\ella+1,\ellb)$, delete the $(q+1)$st $a$-loop (following the order of vertex labels) and root the graph at the vertex where that loop used to be.
\end{itemize}

This algorithm can be modified as in Remark~\ref{rk: drawing cyclically reduced by cb type} to draw uniformly at random a $\PSL$-cyclically reduced subgroup of a given isomorphism type.

\subsection{Implementation and complexity remarks}

\subsubsection*{Two models of computation to measure complexity}

For the complexity analysis, we consider two classical models: the \emph{unit-cost model} (also known as \emph{RAM model}) where each elementary operation, including operations on integers, takes $O(1)$ time; and the \emph{bit-cost} model where an integer $N$ is encoded using $O(\log N)$ space, the number of bits of its representation, and where arithmetic operations are not performed in constant time anymore. This is more realistic in our settings because we are led to handling large integers. For instance, $O(n\log n)$ bits are required to represent the number of size $n$ silhouette graphs (see \cite[Proposition 8.18]{2021:BassinoNicaudWeil}). To simplify the discussion below, we use the following classical notation: for any $\alpha\geq 0$, a non-negative sequence $u_n$ is in $\widetilde O(n^\alpha)$ if there exist constants $C,\beta>0$ such that $u_n\leq C\, n^\alpha \log^\beta n$ for all $n$ sufficiently large. Informally, it means that $u_n$ is in $O(n^\alpha)$ ``up to a poly-logarithmic factor''. Note that, in the bit-cost model, adding or multiplying two numbers encoded with at most $N$ bits costs $\widetilde O(N)$ time. It is elementary for addition and a consequence of, for instance, Harvey and van der Hoeven's result~\cite{2021:HarveyVanDerHoeven} for multiplication. Comparing two numbers encoded with at most $N$ bits costs $O(N)$ time.

In the unit-cost model, we consider that, for any positive integer $n$, we can generate uniformly at random an integer in $[n]$ in $O(1)$ time. In the bit-cost model, we consider that we can generate uniformly at random a bit value of $\{0,1\}$ in $O(1)$ time. If $n$ is a positive integer encoded with $N$ bits, we can therefore produce an element of $[n]$ uniformly at random using a rejection algorithm consisting in repeatedly generating a number made of $N$ independent random bits until the result is in $[n]$. The expected running time of this algorithm is $O(N)$ as the expected number of attempts is at most $2$.
 
\subsubsection*{Precomputing}
Since the parameters $n$, $\ella$, $\ellb$, $\ka$ and $\kb$ are non-negative and satisfy $n = 2\ka+\ella$ and $n \ge 2\kb+\ellb$, there are at most $n^4$ non-zero values for $s(n, \ka,\kb,\ella,\ellb)$ for a given positive integer $n$. They can be computed recursively using Equations~(\ref{eq: fromStb}), (\ref{eq: fromSta}) and~(\ref{eq: fromSy}), the base cases being either trivial (for $n\geq 2$) or given by Proposition~\ref{prop: counting silhouettes}. This yields an $O(n^4)$ time and space algorithm in the unit-cost model and $\widetilde O(n^5)$ time and space algorithm in the bit-cost model.
 
\subsubsection*{Random generation} 
We assume in this section that all the required values of $s(n,\ka,\kb,\ella,\ellb)$ have been precomputed and are accessible in time $O(1)$ in the unit-cost model, $\widetilde O(1)$ in the bit-cost model.

Algorithm $\rcrg$ was written with the proof of Theorem~\ref{thm: proof of algorithm rcrg} in mind. This is the reason why, in particular, it calls for randomly choosing an integer $v$, or integers $v$ and $w$. One can also choose $v = n$, or $v = n$ and $w = n-1$, that is, disregard the randomness of the labeling of the graph we constuct, and add a very last step to the algorithm, which relabels $\Gamma$ by a random permutation. This is a classic trick in the literature on the random generation of labeled combinatorial objects (see for instance~\cite[footnote on p. 12]{1994:FlajoletZimmermannVanCutsem}).
 
In the unit-cost model, $\rsilh$ runs in $O(n)$ average time, using the Fisher-Yates shuffling algorithm~\cite[p.145]{1989:Knuth} and the fact that the number of iterations in $\rsilh$ is bounded in expectation, see Section~\ref{sec: drawing a silhouette graph}. If we use the trick mentioned above (relabeling the graph at the end), every call of the function $\rcrg$ is performed in $O(1)$. As each call decreases the value of $n$ by at least $1$, Algorithm $\rcrg$ runs in $O(n)$ expected time. 

In the bit-cost model, observe that $n$ is encoded using $O(\log n)$ bits, so that all arithmetic operations on $n$, $\ella$, $\ellb$, $\ka$, $\kb$ are performed in $\widetilde O(1)$ time. The bottleneck for the running time of the algorithm lies therefore in
Lines~\ref{rcrg:random x}-\ref{rcrg:compare x}, as generating $x$ and comparing $x$ with the threshold in Line~\ref{rcrg:compare x} both cost $\widetilde O(n)$ time. The overall expected running time of the algorithm in the bit-cost model is therefore $\widetilde O(n^2)$.

This process (generating $x$ and comparing it with a theshold) can be improved using the following idea. Assume that we have two large integers $s$ and $t$, and the sum $s+t$ has already been computed. 
Let $z_0\cdots z_{N-1}$ be the binary encoding of $s+t$ (with $z_0 = 1$). Let also $s_0\cdots s_{N-1}$ be the binary encoding of $s$ (here $s_0$ may be 0). Generating $x$ in $[s+t]$ and comparing it to $s$, amounts to simulating a Bernoulli law of parameter $\frac{s}{s+t}$ (in the bit-cost model). This is performed using Algorithm~$\berna$, which generates a uniform random integer, say $x$, between $0$ and $2^N-1$ bit by bit, halting as soon as we are guaranteed that one of the three possible situations holds : $x \geq s+t$ (Failure), $x<s$ (\texttt{True}) and $s\leq x<s+t$ (\texttt{False}).

\begin{algorithm}[H]
\DontPrintSemicolon
$smaller_{s+t} = \emptyset$\;
$smaller_{s} = \emptyset$\;
\For{$i\in\{0,\ldots, N-1\}$}{
	$bit$ = Uniform($\{0,1\}$)\;
	\If{$smaller_{s+t} = \emptyset$\label{berna:x+y}}
	{ 
		\lIf{$bit>z_i$}{\Return Failure}
		\lIf{$bit<z_i$}{$smaller_{s+t} = \texttt{True}$}
	}
	\If{$smaller_{s} = \emptyset$\label{berna:x}}
	{ 
		\lIf{$bit>s_i$}{$smaller_{s} = \texttt{False}$}
		\lIf{$bit<s_i$}{$smaller_{s} = \texttt{True}$}
	}
	\If{$smaller_{s+t} \neq \emptyset$ and $smaller_{s} \neq \emptyset$}{\Return $smaller_{s}$}
}
\If{$smaller_{s+t} = \emptyset$}{\Return Failure \tcp{the generated number is  $\geq s+t$}}
\Return  \texttt{False} \tcp{the generated number is equal to $s$}
\caption{$\berna(x,y,N)$}\label{algo:berna}
\end{algorithm}
\medskip

The main algorithm to simulate the Bernoulli law of parameter $\frac{s}{s+t}$ consists in repeatedly calling $\berna$ until the result is not Failure. The analysis of the expected number of bits generated in the process is straightforward, as $smaller_{s+t}$ and $smaller_s$ are determined with probability $\frac12$ at each iteration of the first loop: the number of iterations required to determine each one of them is bounded above by a geometric law of parameter $\frac12$. Since $s+t > 2^{N-1}$, the expected number of calls to $\berna$ is bounded above by a constant. Hence the expected bit-cost complexity of our procedure to simulate the Bernoulli law is $\widetilde\O(1)$.

To implement the announced improvement, we modify the precomputation step by storing not only the values of $s(n,\ka,\kb,\ella,\ellb)$, but also their bit-lengths and the values $n\cdot (\kb+1)\cdot s(n-1,\ka,\kb+1,\ella-1,0)$ (used Line~\ref{rcrg:compare x}). Simulating a Bernoulli law of parameter $n\cdot (\kb+1)\cdot s(n-1,\ka,\kb+1,\ella-1,0)/s(n,\ka,\kb,\ella,\ellb)$ instead of performing Lines~\ref{rcrg:random x} and~\ref{rcrg:inner if}, lowers the expected time complexity of $\rcrg$ to $\widetilde O(n)$.
 
\begin{remark}
When $\ellb>0$, one can directly choose  $\ell_3$ $a$-loops of the graph $\Delta$ built at Line~\ref{rcrg:delta ellb} to apply the inverse of a $\lambda_3$-move $\ell_3$ times directly. This does not change the overall complexity of $\rcrg$ in both models of computation. Similarly, if $\ella=\ellb=0$ and $\kb=n/2$, the generated graph is a cycle made of an alternation of $a$-transitions and $b$-transitions, which could be generated directly without making several recursive call. Again, this does not change the complexity.
\end{remark}


{\small\bibliographystyle{abbrv}

}

\end{document}